\documentstyle[amscd,amssymb,verbatim]{amsart}

\makeatletter
\def\@maketitle{%
\defaultfont\normalsize
\let\@makefnmark\relax \let\@thefnmark\relax \ifx\@empty\@subjclass\else
\@footnotetext{1991 {\it Mathematics Subject
Classification}.\enspace
\@subjclass.}\fi
\ifx\@empty\@keywords\else
\@footnotetext{{\it Key words and phrases.}\enspace \@keywords.}\fi
\ifx\@empty\@thanks\else
\@footnotetext{\@thanks}\fi
\topskip66\p@ 
\vtop{\centering{\baselineskip14\p@\bf
\expandafter{\@title}\@@par}%
\global\dimen@i\prevdepth}%
\prevdepth\dimen@i
\ifx\@empty\@authors
\else
\baselineskip32\p@
\vtop{\@andify{ AND }\@authors
\centering{{\@authors}\@@par}%
\global\dimen@i\prevdepth}\relax
\prevdepth\dimen@i
\fi
\ifx\@empty\@dedicatory
\else
\baselineskip18\p@
\vtop{\centering{\small\it\@dedicatory\@@par}%
\global\dimen@i\prevdepth}\prevdepth\dimen@i \fi
\ifx\@empty\@date\else
\baselineskip24\p@
\vtop{\centering\@date\@@par
\global\dimen@i\prevdepth}\prevdepth\dimen@i \fi
\normalsize
\dimen@32\p@ \advance\dimen@-\baselineskip \vskip\dimen@\@plus14\p@
} 
\makeatother

\theoremstyle{plain}
\newtheorem{Thm}{Theorem}[section]
\newtheorem{Cor}{Corollary}[section]
\newtheorem{Lem}{Lemma}[section]
\newtheorem{Prop}{Proposition}[section]

\theoremstyle{definition}
\newtheorem{Def}{Definition}[section]
\newtheorem{Ex}{Example}[section]
\theoremstyle{remark}

\newtheorem{Rems}{Remarks}[section]

\renewcommand{\a}{\alpha}
\renewcommand{\b}{\beta}
\newcommand{\ra}{\rightarrow}
\newcommand{\Spec}{\operatorname{Spec}}
\newcommand{\sub}{\subset}
\newcommand{\wt}{\widetilde}
\renewcommand{\O}{{\cal O}}
\newcommand{\NN}{{\cal N}}
\newcommand{\Der}{\operatorname{Der}}
\newcommand{\Hom}{\operatorname{Hom}}
\newcommand{\NCS}{\operatorname{NCS}}
\newcommand{\lrar}[1]{\begin{picture}(50,10)(-25,-5)
\put(-25,0){\vector(1,0){50}}
\put(0,5){\makebox(0,0)[b]{\mbox{$#1$}}} \end{picture}}

\newcommand{\ldar}[1]{\begin{picture}(10,50)(-5,-25)
\put(0,25){\vector(0,-1){50}}
\put(5,0){\mbox{$#1$}}
\end{picture}}

\newcommand{\lbrurar}[1]{\begin{picture}(50,50)(-25,-25)
\multiput(-25,-25)(5,5){10}{\circle*{1}} \put(25,25){\vector(1,1){0}}
\put(5,-5){\mbox{$#1$}}
\end{picture}}

\def\morph#1{\overset{#1}{\ra}}
\def\rings{\text{\bf Rings}}
\def\dsp#1{$\displaystyle{#1}$}

\def\ts{\times_S}
\def\boxs{\boxtimes_S}
\def\derboxs{{\overset{\Bbb L}{\boxtimes}}_S} \def\ots{\otimes_S}

\theoremstyle{remark}
\newtheorem{notation}{Notation}[section] 

\numberwithin{equation}{subsection}

\renewcommand{\O}{{\cal O}}
\newcommand{\D}{{\cal D}}

\newcommand{\Sym}{\operatorname{Sym}}
\newcommand{\hra}{\hookrightarrow}

\newcommand{\End}{\operatorname{End}}
\newcommand{\Ext}{\operatorname{Ext}}
\newcommand{\Proj}{\operatorname{Proj}}
\newcommand{\Nm}{\operatorname{Nm}}
\newcommand{\rk}{\operatorname{rk}}
\newcommand{\id}{\operatorname{id}}
\newcommand{\sym}{\operatorname{sym}}
\newcommand{\gr}{\operatorname{gr}}
\newcommand{\ml}{\operatorname{ml}}
\newcommand{\PP}{\Bbb P}
\newcommand{\C}{\Bbb C}
\newcommand{\Z}{\Bbb Z}

\title{Fourier transform for D-algebras} 
\author{A. Polishchuk and M. Rothstein}
\thanks{This research was supported in part by
the NSF grants DMS-9700458 and  DMS-9626522}

\begin{document}

\maketitle

\bigskip

This paper is devoted to the construction of an analogue of the Fourier
transform for a certain
class of non-commutative algebras. The model example which initiated this
study is the equivalence
between derived categories of $\D$-modules on an abelian variety and
$\O$-modules on the universal
extension of the dual abelian variety by a vector space (see [L], [R2]).
The natural framework for
a generalization of this equivalence is provided by the language of
$D$-algebras developed by
A.~Beilinson and J.~Bernstein in [BB]. We consider a subclass of
$D$-algebras we call special $D$-algebras.
We show that whenever
one has an equivalence of categories of $\O$-modules on two varieties $X$
and $Y$, it gives rise to a correspondence between special
$D$-algebras on $X$ and $Y$ such
that the corresponding derived categories of modules are equivalent. When
$X$ is an abelian
variety, $Y$ is the dual abelian variety, according to Mukai [M] the
categories of $\O$-modules on
$X$ and $Y$ are equivalent, so our construction gives in particular the
Fourier transform between
modules over rings of twisted differential operators ({\it tdo}, for short,
see section
\ref{gen-sec} for a definition) with non-degenerate first Chern class on
$X$ and $Y$.

We also deal with the microlocal version of the Fourier transform. The
microlocalization
of a special filtered $D$-algebra on $X$ is an NC-scheme in the sense of
Kapranov (see [K]) over
$X$, i.e. a ringed space whose structure ring is complete with respect to
the topology defined by
commutator filtration (see section \ref{et-sec}). We show that in our
situation the derived
categories of coherent sheaves on microlocalizations are also equivalent.
In the case of rings of
twisted differential operators on the dual abelian varieties $X$ and $Y$
one can think about the
corresponding microlocalized algebras as deformation quantizations of the
cotangent spaces $T^*X$
and $T^*Y$. The projections $p_X:T^*X\ra T_0^*X$ and $p_Y:T^*Y\ra T_0^*Y$
can be considered as
completely integrable systems with dual fibers (a choice of a
non-degenerate tdo on $X$ induces an
identification of the bases $T_0^*X$ and $T_0^*Y$ of these systems). We
conjecture that our
construction generalizes to other dual completely integrable systems. This
hope is based on the
following observation: the relative version of our transform gives a
Fourier transform for modules
over relative tdo's on dual families of abelian varieties, while
deformation quantizations are
usually subalgebras in these tdo's.
An example of dual completely integrable systems appears in the geometric
Langlands program (see [BD]). Namely, one may consider Hitchin systems for
Langlands dual
groups. The analogue of the Fourier transform in this situation should lead
to the equivalence
between modules over a microlocalized tdo on moduli spaces of principal
bundles for Langlands
dual groups.

An important aspect of our work is that in an appropriate sense, the
microlocalized Fourier
transform is in fact \'etale local.	Given a special filtered $D$-algebra
$\cal A$ on $X$, let $\cal A_{ml}$ be the corresponding microlocalized
scheme. Then $\cal A_{ml}$
is a non-commutative thickening of the product of $X$ with a scheme $Z$.
Denoting by $\Phi(\cal A)$ the corresponding $D$-algebra on $Y$, $\Phi(\cal
A)_{ml}$
is a thickening of $Y\times Z$. We then prove that the microlocalized Fourier
transform is \'etale local in $Z$. It should be noted that $\cal A_{ml}$
and $\Phi(\cal A)_{ml}$ are not schemes over $Z$, so one does not have a
straightforward base-change
argument. Rather,
%
%
we develop in section
\ref{et-sec} the non-commutative version of the theory of \'etale morphisms
in the framework of
Kapranov's NC-schemes, and establish the equivalence by a version of the
topological invariance of
\'etale morphisms.

Our work is motivated in part by Krichever's construction of solutions of
the KP-hierarchy, [Kr]. (See section \ref{geom-subsec}.)
%
Let $W$ be a smooth variety
of dimension $r$, embedded in its Albanese variety, $X$. Let
$D\sub W$ be an ample hypersurface and
$V\sub H^0(D,\O_D(D))$ an $r$-dimensional basepoint free subspace. Let
$\phi:D\ra\PP(V^*)$ denote the corresponding morphism and let $U\sub D$ be
an open subset such that $\phi|_U$ is \'etale. Let $Y=Pic^0(W)$. Then
$V$ maps to the space of vector fields on $Y$, and hence one has a
subalgebra of the differential
operators on $Y$, consisting of those operators which differentiate only in
the ``$V$" directions.
Denote this algebra by $\Phi(\cal A_V)$. That is, $\Phi(\cal A_V)$ is dual to a
$D$-algebra $\cal A_V$ on $X$. Then the microlocalizations of these
$D$-algebras are
thickenings of $Y\times\PP(V^*)$ and $X\times\PP(V^*)$ respectively. In
particular,
one has the
\'etale localizations $\cal A_{ml,U}$ and $\Phi(\cal A)_{ml,U}$ supported
on $X\times U$ and
$Y\times U$ respectively. Let $U_{\infty}$ denote the formal neighborhood
of $U$ in $W$. Then
the diagonal embedding $U\rightarrow X\times U$ extends to an embedding
$\Delta_{\infty}:U_{\infty}\rightarrow
\cal A_{ml,U}$. On the other hand, $U\times Y$ sits in both $X\times Y$ and
$\Phi(\cal A)_{ml,U}$. Denote by $\cal L_{\infty}$ the Fourier transform of
$\Delta_{\infty *}(\cal O_{U_{\infty}})$. We prove that $\cal L_{\infty}$
is a locally free
rank-one left $\cal O$ module on $\Phi(\cal A)_{ml,U}$, whose restriction
to $U\times Y$ is the
restriction of the Poincar\'e line bundle. Thus $\cal L_{\infty}$ is a
deformation
of the Poincar\'e line bundle. Furthermore, for any positive integer $k$,
	the Fourier
transform of $\Delta_{\infty *}(\cal O_{U_{\infty}}(kU))$ is $\cal
L_{\infty}(k(U\times Y))$.

The point is that the ring $H^0(W,\cal O(*D))$ acts by $\cal
A_{ml,U}$-endomorphisms on \hfill\break
$\Delta_{\infty *}(\cal O_{U_{\infty}}(*U))$. Functoriality of the Fourier
transform then
gives us a representation
\begin{equation}\label{representation}
H^0(W,\cal O(*D))\rightarrow End_{\Phi(\cal A)_{ml,U}}(\cal
L_{\infty}(*(U\times Y)))\ .\end{equation}
When $W$ is a curve and $D$ is a point, this representation reduces to the
Burchnall-Chaundy [BC] representation of $H^0(W,\cal O(*D))$ by
differential operators. We intend to study
the representation (\ref{representation}) further in a future work. In
particular, the problem of
characterizing the image of this representation is quite interesting, and
should lead to
generalizations of the KP-hierarchy.

\bigskip

\begin{notation}
Fix a scheme $S$. Given an $S$-scheme $U$, denote by $\pi^U_S$ the structural
morphism. By ``associative $S$-algebra on $U$" we mean a sheaf of
associative rings $\cal A$ on
$U$ equipped with a morphism of sheaves of rings from ${\pi^U_S}^{-1}(\cal
O_S)$ to the center of $\cal A$. We abbreviate
``$\otimes_{{\pi^U_S}^{-1}(\cal O_S)}$" by ``$\ots$". For a scheme $U$
we denote by
$\cal D^b(U)$ the bounded derived category of quasicoherent sheaves on $U$.
Throughout the
paper, $X$ and $Y$ are flat, separated $S$-schemes. \end{notation}

\section{$D$-algebras and Lie algebroids}\label{gen-sec} \subsection{}
Let us recall some definitions from [BB]. A {\it differential} $\cal
O_X$-{\it bimodule} $M$ is a quasicoherent sheaf on $X \ts X$ supported on
the diagonal $X \subset X \ts X$. One can consider the category of
differential $\cal O_X$-bimodules as a subcategory in the category of all
sheaves of $\cal O_X$-bimodules on $X$. A $D$-{\it algebra} on $X$ is a
sheaf of flat, associative $S$-algebras $\cal A$ on $X$ equipped with a
morphism of $S$-algebras $i: \cal O_X \rightarrow \cal A$ such that $\cal
A$ is a differential $\cal O_X$ - bimodule. This means that $\cal A$ has an
increasing filtration $0=\cal A_{-1} \subset \cal A_0 \subset \cal A_1
\subset \cdots$ such that $\cal A =
\cup \cal A_n$ and $ ad(f) (\cal A_k) \subset \cal A _{k-1}$ for any $k \ge
0$ and $f \in \cal O_X$ where $ ad(m): = rm-mr$. We denote by $b(\cal A)$
the quasi-coherent sheaf on $X\ts X$ (supported on the diagonal)
corresponding to $\cal A$.
Also we denote
by $ \cal M ( \cal A)$ the category of sheaves of left $\cal A$-modules on
$X$ which
are quasicoherent as $\cal O_X$-modules.

\subsection{}\label{circle} Let us describe some basic operations with
$D$-algebras and modules over them. Let
$\cal A_X$ and $\cal A_Y$ be $D$-algebras over $X$ and $Y$ respectively.
One defines a $D$-algebra
$\cal A_X\boxs
\cal A_Y$
on $X\ts Y$ by gluing $D$-algebras over products of affine opens $U\ts V$
corresponding to $\cal A_X(U)\ots \cal A_Y(V)$. A module $M\in\cal M(\cal
A_X\boxs \cal A_Y)$ is the same as a quasicoherent $\cal O_{X\ts Y}$-module
together with commuting actions of $p_X^{-1}(\cal A_X)$ and $p_Y^{-1}(\cal
A_Y)$
which are compatible with the $\cal O_{X\ts Y}$-module structure (where
$p_X$ and $p_Y$ are projections from $X\ts Y$ to $X$ and $Y$) . In
particular, we have the natural structure of $D$-algebra on
$p_X^*\cal A_X\simeq \cal A_X\boxs\cal O_Y$ and $p_Y^*\cal A_Y\simeq \cal
O_X\boxs\cal A_Y$ and natural embeddings of $D$-algebras $p_X^*\cal
A_X\hookrightarrow \cal A_X\boxs \cal A_Y$, $p_Y^*\cal A_Y\hookrightarrow
\cal A_X \boxs\cal A_Y$.
For a pair of modules $M_X\in\cal M(\cal A_X)$ and $M_Y\in\cal M(\cal A_Y)$
there is a natural structure of $\cal A_X\boxs\cal A_Y$-module on $M_X\boxs
M_Y$.

Now assume that we have $D$-algebras $\cal A_X$, $\cal A_Y$, and $\cal A_Z$
on $X$, $Y$ and $Z$ respectively. Then we can define an operation
$$\circ_{\cal A_Y}:\cal D^-(\cal M(\cal A_X\boxs\cal A_Y^{op}))\times \cal
D^-(\cal M(\cal A_Y\boxs\cal A_Z))\rightarrow \cal D^-(\cal M(\cal
A_X\boxs\cal A_Z)).$$ The definition is the globalization of the operation
of tensor product of bimodules.
Namely, for a pair of objects $M\in \D^-(\cal M(\cal A_X\boxs \cal
A_Y^{op}))$ and $N\in \D^-(\cal M(\cal A_Y\boxs \cal A_Z))$ we can form the
external tensor product $M\boxs N\in
\D^-(\cal M(\cal A_{XYZ}))$ where $\cal A_{XYZ}= \cal A_X\boxs \cal
A_Y^{op}\boxs\cal A_Y\boxs\cal A_Z$ is a $D$-algebra on $X\ts Y\ts Y\ts Z$.
Note that there is a natural structure of left $\cal A_Y\boxs\cal
A_Y^{op}$-module on
$b(\cal A_Y)$ given by the multiplication in $\cal A_Y$. Hence, we can
consider the tensor product $$(\cal A_X\boxs
b(\cal A_Y)\boxs\cal A_Z) \overset{\Bbb L}{\otimes}_ {\cal A_{XYZ}}
	(M\derboxs N)$$
as an object in the category
$\D^-(p_{XZ}^{-1}(\cal A_X\boxs\cal A_Z))$ where $p_{XZ}^{-1}$ denotes a
sheaf-theoretical inverse image. 
Finally, we set
$$M\circ_{\cal A_Y}N=Rp_{XZ*}((\cal A_X\boxs b(\cal A_Y)\boxs\cal A_Z)
\overset{\Bbb L}{\otimes}_
{\cal A_{XYZ}} (M\derboxs N))
.$$

There is also the following equivalent definition: $$M\circ_{\cal A_Y}N=
Rp_{XZ*}((M\boxs\cal A_Z)\overset{\Bbb L}{\otimes}_ {{\cal
A_X}^{op}\boxs\cal A_Y\boxs\cal A_Z}({\cal A_X}^{op}\boxs N))\ .$$



Specializing to the case that $Z=S$ and $\cal A_Z=\cal O_S$, we see that
every $\cal A_X\boxs\cal A_Y^{op}$-module $F$ on $X\ts Y$ defines a functor
$G\mapsto F\circ_{\cal A_Y} G$
from $\cal D^-\cal M(\cal A_Y)$ to $\cal D^-\cal M(\cal A_X)$.

\begin{Prop}
\begin{enumerate}
\item The operation $\circ$ is associative in the natural sense.
\vskip 5pt
\item One has $M\circ b(\cal A_Y)\simeq M$ and $b(\cal A_Y)\circ N\simeq N$
canonically, for \hfill\break\hbox{$M\in \cal D^-\cal M(\cal
A_X\boxs\cal A_Y^{op})$} and $N\in \cal D^-\cal M(\cal A_Y\boxs\cal A_Z)$.
\vskip 5pt
\item If $M$ is a differential $\cal O_X$-bimodule and $N$ is an $\cal
O_{X\ts Y}$-module,
then\hfill\break \hbox{$
b(M)\circ_{\cal O_X}N=p_X^{-1}(M)\otimes_{p_X^{-1}(\cal O_X)}N$.}
\end{enumerate}
\end{Prop}

It follows that if $\cal A$ is a $D$-algebra on $X$, the structural morphism
on $\cal A$ may be
viewed as a morphism
$b(\cal A)\circ_{\cal O_X} b(\cal A)\rightarrow b(\cal A)$, and if $M$ is a
left $\cal A$-module,
the action of $\cal A$ on $M$ is given by a morphism $b(\cal A)\circ_{\cal
O_X} M\rightarrow M$.
Moreover, an
$\cal A_X\boxs\cal A_Y^{op}$-module structure on an $\cal O_{X\ts
Y}$-module $M$ is the same
as a pair of morphisms $b(\cal A_X)\circ_{\cal O_X} M\rightarrow M$ and
$M\circ_{\cal O_Y}b(\cal
A_Y)\rightarrow M$ making
$M$ a (left $b(\cal A_X)$)-(right $b(\cal A_Y)$)-module with respect to
$\circ$, such that the
two module structures commute.

\subsection{}
Recall that a {\it Lie algebroid} $L$ on $X$ is a (quasicoherent) $\cal
O_X$-module equipped with a morphism of $\cal O_X$-modules $\sigma: L
\rightarrow \cal T (:= \text{Der}_S \cal O_X = $\text{relative tangent
sheaf of $X$}) and an $S$-linear Lie bracket
$[\cdot,\cdot] : L\ots L \rightarrow L$ such that $\sigma$ is a
homomorphism of Lie algebras and the following identity is satisfied:
$$
[\ell_1, f \ell_2] = f\cdot [\ell_1, \ell_2] + \sigma (\ell_1) (f) \ell_2
$$
where $\ell_1, \ell_2 \in L, f \in \cal O_X$. To every Lie algebroid $L$
one can associate a $D$-algebra $\cal U (L)$ called the {\it universal
enveloping algebra} of $L$. By definition $\cal U(L)$ is a sheaf of
algebras equipped with the morphisms of sheaves $i: \cal O_X \rightarrow
\cal U(L) , i_L: L \rightarrow \cal U(L)$, such that $\cal U(L)$ is
generated, as an algebra, by the images of these morphisms and the only
relations are:

\begin{enumerate}
\renewcommand{\labelenumi}{(\roman{enumi})} \item $i$ is a morphism of
algebras;
\item $i_L$ is a morphism of Lie algebras; \item $i_L(f \ell)= i(f) i_L
(\ell),\ [ i_L (\ell) , i(f)] = i (\sigma (\ell) (f))$, where $f \in \cal
O_X, \ell \in L$. \renewcommand{\labelenumi}{(\arabic{enumi})}
\end{enumerate}

\subsection{}
Let $L$ be a Lie algebroid on $X$. A central extension of $L$ by $\cal O_X$
is a Lie algebroid $\tilde{L}$ on $X$ equipped with an embedding of $\cal
O_X$ - modules $c:\cal O_X \hookrightarrow \tilde{L}$ such that $[c(1),
\tilde{\ell}] = 0$ for every $\tilde{\ell} \in \tilde{L}$ ( in particular,
$c(\cal O_X)$ is an ideal in $\tilde{L}$), and an isomorphism of Lie
algebroids $\tilde{L}/c(\cal O_X) \simeq L$. For such a central extension
we denote by $\cal U^{\circ}(\tilde{L})$ the quotient of $\cal
U(\tilde{L})$ modulo the ideal generated by the central element $i(1) -
i_{\tilde{L}} (c(1))$.

\begin{Lem}\label{locallyfree}
Let $L$ be a locally free $\cal O_X$-module of finite rank. Then there is a
bijective correspondence between isomorphism classes of the following data:
\begin{enumerate}
\renewcommand{\labelenumi}{(\roman{enumi})} \item a structure of a Lie
algebroid on $L$ and a central extension $\tilde{L}$ of $L$ by $\cal O_X$.
\item a $D$-algebra $\cal A$ equipped with an increasing algebra filtration
$\cal O_X = \cal A_0 \subset \cal A_1 \subset \cal A_2 \subset \dots$ such
that $\cup \cal A_n = \cal A$ and an isomorphism of the associated graded
algebra $\gr A$ with the symmetric algebra $S^\bullet L$. \end{enumerate}
\renewcommand{\labelenumi}{(\arabic{enumi})} \end{Lem}
The correspondence between (i) and (ii) maps a central extension $
\tilde{L}$ to $\cal U^\circ (\tilde{L})$.

\subsection{}
Assume that $X$ is smooth over $S$. Then one can take $L= \cal T$ with its
natural Lie algebroid structure. The corresponding central extensions
$\tilde{\cal T}$ of $\cal T$ by $\cal O$ are called {\it Picard algebroids}
and the associated $D$-algebras are called {\it algebras of twisted
differential operators}; or simply {\it tdo}'s. If $\cal D$ is a tdo, $\cal
D_{-1} = 0 = \cal D_{0} \subset \cal D_1 \subset \cal D_2 \subset
\dots$ is its maximal $D$-filtration, i.e. $$
\D_i = \{d \in \D | ad(f) d \in \D_{i-1}\ ,\ f\in \cal O_X \},$$ then $\gr\D
\simeq S^\bullet \cal T$.
\begin{Lem}\label{isomorphism}
For a locally free $\cal O_X$-module of finite rank $E$ one has a canonical
isomorphism
$$
\Ext^1_{\cal O_{X \ts X}}(\Delta_* E, \Delta_*\cal O_X) \simeq \Hom_{\cal
O_X} (E, \cal T) \oplus \Ext^1_{\cal O_X} (E,\cal O_X), $$
where $X\overset \Delta\rightarrow X\ts X$ is the diagonal embedding. \end{Lem}

\begin{pf}
Since $\Delta_* E \simeq p^*_1 E \otimes_{\cal O_{X\ts X}} (\cal O_{X\ts
X}/J)$,
where $J$ is the ideal sheaf of the diagonal, we have an exact sequence
$$
0 \rightarrow \Hom (p^*_1 E \otimes_{\cal O_{X\ts X}} J,\Delta_*\cal O_X)
\rightarrow \Ext^1 (\Delta_*E, \Delta_*\cal O_X) \rightarrow \Ext^1 (p^*_1
E,\Delta_*\cal O_X)
$$
Note that the first and last terms are isomorphic to $\Hom(E, \cal T)$ and
$\Ext^1(E,\cal O_X)$ respectively. It remains to note that there is a
canonical splitting $\Delta_*: \Ext^1(E,\cal O_X) \rightarrow
\Ext^1(\Delta_*
E, \Delta_* \cal O)$.
\end{pf}

Note that the projection $\Ext^1(\Delta_*E, \Delta_*\cal O_X) \rightarrow
\Hom(E, \cal T)$ can be described as follows. Given an extension
\[
\begin{CD}
0 @>>> \Delta_* \cal O_X @>>> \tilde{E} @>>> \Delta_* E @>>> 0 \end{CD}
\]
the action of $J/J^2$ on $\tilde{E}$ induces the morphism $J/J^2 \otimes
\tilde{E} \rightarrow \Delta_* \cal O$, which factors through $J/J^2
\otimes \Delta_* E$, since $J$ annihilates $\Delta_* \cal O$. Hence we get
a morphism $\Delta_* E \rightarrow \Delta_* \cal T$.

Now if $\cal A$ is a $D$-algebra,
equipped with a filtration $\cal
A_{\bullet}$ such that
$\gr\cal A \simeq S^{\bullet}(E)$, then we consider the corresponding
extension of $\cal O_X$-bimodules \[
\begin{CD}
0 @>>> \cal O_X = \cal A_0 @>>> \cal A_1 @>>> E=\cal A_1/\cal A_0@>>> 0
\end{CD}
\]
as an element in $\Ext^1_{\cal O_{X \times X}}(\Delta_* E, \Delta_* \cal
O_X)$. By definition, $\cal A$ is a tdo if the projection of this element
to $\Hom_{\cal O_X}(E, \cal T)$ is a map $E \rightarrow \cal T$ which is an
isomorphism.

\section{Equivalences of categories of modules over $D$-algebras}

\subsection{} Let
$P$ be
an object in $\D^b(X\ts Y)$, $Q$ be an object in $\D^b(Y\ts X)$ such that
$$P\circ_{\cal O_Y} Q\simeq\Delta_*\cal O_Y,\\ Q\circ_{\cal O_X}
P\simeq\Delta_*\cal O_X.$$ where $\Delta$ denotes the diagonal
embedding. In this case the functors
$$\Phi_P:M\mapsto P\circ M,\\
\Phi_Q:N\mapsto Q\circ N$$ establish an equivalence of categories $\D^-(X)$
and $\D^-(Y)$.

For example, we have these data in the following situation: $X$ is an
abelian $S$-scheme,
$Y=\hat{X}$ is the dual abelian $S$-scheme, $P=\cal P$ is the normalized
Poincar\'e line bundle on $X \ts \hat{X}$, $Q=\sigma^*\cal
P^{-1}\omega_{X}^{-1}[-g]$ where $\sigma:\hat{X}\ts X\ra X\ts \hat{X}$ is
the permutation of factors, $g = \dim X$.

\subsection{} Let us call a quasicoherent sheaf $K$ on $X \ts X$ {\it
special} if there is a
filtration $0 = K_{-1}\subset K_0 \subset K_1 \subset \dots$ of $K$ and a
sequence of sheaves
of flat, quasicoherent $\cal O_S$-modules $F_i$ such that $\cup K_i = K$
and $K_i/K_{i-1} \simeq \Delta_* {\pi_S^X}^*(F_i)$ for every $i \ge 0$. We
denote by $\cal S_X$ the exact category of special sheaves on $X \ts X$.
The following properties are easily verified. \begin{Lem}\label{sheaves on
S} Let $F\in \D^-(\cal M(\cal O_S))$. \begin{enumerate}
\item Let $G\in D^-(\cal M(\cal O_{X\ts Y}))$. Then
$\Delta_{X*}{\pi_S^X}^*(F)\circ_{\cal O_X}G= {\pi_S^{X\ts
Y}}^*(F)\overset{\Bbb L}{\otimes}_{\cal O_{X\ts Y}}G$. \item $P
\circ_{\cal O_Y} {\pi_S^Y}^*(F) \circ_{\cal O_Y} Q={\pi_S^X}^*(F)$.
\end{enumerate}
\end{Lem}
The following proposition then follows.

\begin{Prop}
\begin{enumerate}
\item
For every $K \in \cal S_Y$, the functor $\cal M(\cal O_{Y\times
Z})\rightarrow \cal M(\cal O_{Y\times Z})$, $M\mapsto K\circ_{\cal O_Y} M$,
is exact. \vskip 5pt
\item
For every
pair of special sheaves $K, K' \in \cal S_Y$, $ K{\circ}_{\cal O_Y} K'$ is
special.
\end{enumerate}
\end{Prop}

\begin{Prop} The functor $\Phi: K \mapsto P \circ_{\cal O_Y} K \circ_{\cal
O_Y} Q$
defines an equivalence of categories $\Phi: \cal S_Y \rightarrow \cal S_X$.
\end{Prop}
\begin{pf}
>From lemma \ref{sheaves on S}, together with the fact that the
operation $\circ$ commutes with inductive limits, we have $\Phi(K) \in \cal
S_X$ for every $K \in
\cal S_Y$. It remains to notice that there is an inverse functor to $\Phi$
given by $$
\Phi^{-1}(K')=Q {\circ}_{\cal O_X} K' {\circ}_{\cal O_X} P $$
where $K' \in \cal S_X$.
\end{pf}

\begin{Prop}\label{homomorphism}
For $K, K' \in \cal S_Y, M \in \cal D^b(Y)$ one has a canonical isomorphism
of $\cal O_X$-bimodules $$
\Phi(K \circ_{\cal O_{Y}} K') \simeq \Phi K \circ_{\cal O_X} \Phi K', $$
and a canonical isomorphism in $\cal D^b(X)$ $$ \Phi (K \circ_{\cal O_Y} M)
\simeq \Phi K \circ_{\cal O_X} \Phi M. $$
\end{Prop}

\begin{Def}A $\it special$ $D$-$\it algebra$ on $X$ is a $D$-algebra $\cal
A$ such that
the sheaf $b(\cal A)$ on
$X\ts X$ is special.\end{Def}

It follows from the above proposition that for any special $D$-algebra $
\cal A$ on $Y$
there exists a
canonical $D$-algebra $\Phi \cal A$ on $X$ such that $$b(\Phi\cal A)\simeq
\Phi(b(\cal A)).$$
Namely one just has to apply
$\Phi$ to structural morphisms $b(\cal A) \circ_{\cal O_{Y}} b(\cal A)
\rightarrow b(\cal A)$
and $\Delta_* \cal O_{Y} \rightarrow b(\cal A)$. Futhermore, we now prove
that the derived
categories of modules over $\cal A$ and $\Phi \cal A$ are equivalent.

\begin{Thm}\label{main1} Assume that $P$ and $Q$ are quasi-coherent sheaves
up to a shift (i.e. they have only one cohomology). Then for every special
$D$-algebra $\cal A$ on $Y$ there is a canonical exact equivalence $\Phi:
\cal D^- \cal M(\cal A) \rightarrow \cal D^- \cal M(\Phi \cal A)$ such that
the following diagram of functors is commutative: \end{Thm} \[
\begin{CD}
\cal D^b \cal M(\cal A) @>\Phi>> \cal D^b \cal M(\Phi \cal A)\\ @VVV@VVV\\
\cal D^b (Y) @>\Phi>> \cal D^b (X)
\end{CD}
\]
where the vertical arrows are the forgetting functors. \begin{pf} Let us
consider the following
object in $\D^b(X\ts Y)$: $$ \cal B = P {\circ}_{\cal O_{Y}} b(\cal A) = P
\otimes_{\cal O_{X\ts Y}}
p_{Y}^*\cal A\ . $$ Note that $\cal B$ is actually concentrated in one
degree so we can consider
it as a quasicoherent sheaf on
$X\ts Y$ (perhaps shifted).
We claim that there is a canonical $\Phi \cal A\boxtimes \cal
A^{op}$-module structure on $\cal B$. Indeed, it suffices to construct
commuting actions $b(\Phi \cal A)\circ\cal B\rightarrow\cal B$ and $\cal
B\circ b(\cal A)\rightarrow\cal B$
compatible with $\cal O_{X\ts Y}$-module structure. The right action of
$\cal A$ is obvious
while the left action of $\Phi \cal A$ on $\cal B$ is given by the
following map:
$$
\Phi b(\cal A) \circ_{\cal O_{X}} \cal B = P \circ_{\cal O_{Y}} b(\cal A)
{\circ}_{\cal O_{Y}} Q {\circ}_{\cal O_X} P \circ_{\cal O_{Y}} b(\cal A)
\rightarrow P \circ_{\cal O_{Y}} b(\cal A) \circ_{\cal O_{Y}} b(\cal A)
\rightarrow P \circ_{\cal O_{Y}} b(\cal A) = \cal B\ , $$ where the
last arrow is induced by multiplication in $\cal A$. It is clear that this
is a morphism of right
$p_{Y}^{-1}\cal A$ modules. On the other hand, there is a natural
isomorphism of sheaves on $X\ts
Y$ $$ b(\Phi \cal A)\circ_{\cal O_{X}} P\simeq P\circ_{\cal O_{Y}} b(\cal
A)\circ_{\cal O_{Y}} Q
\circ_{\cal O_X} P\simeq P\circ_{\cal O_{Y}} b(\cal A)\simeq\cal B\ . $$
One can easily check that
the above left action of $p_X^{-1}\Phi \cal A$ on $\cal B$ is compatible
with the natural left
$p_X^{-1}\Phi \cal A$-module structure on $b(\Phi \cal A)\circ_{\cal O_{X}}
P$ via this
isomorphism. Thus, $\cal B$ is an object of $\D^b(\cal M(\Phi\cal
A\boxtimes\cal A^{op}))$
(concentrated in one degree).
So we can define the functor
$$
\Phi: \D^b\cal M(\cal A) \rightarrow \D^b
\cal M(\Phi \cal A) : M \mapsto \cal B{\circ}_{\cal A} M $$
Similarly,
we define an $\cal A\boxtimes\Phi \cal A^{op}$-module (perhaps shifted) on
$Y\ts X$:
$$
\cal B' = b(\cal A) \circ_{\cal O_Y} Q \simeq \cal Q \circ_{\cal O_X}
b(\Phi \cal A)
$$
and the functor
$$
\Phi^\prime: \D^b \cal M(\Phi \cal A) \rightarrow \D^b \cal M(\cal A) : N
\mapsto B^\prime
{\circ}_{\Phi \cal A} N.
$$
One has an isomorphism in the derived
category of right $\cal O_Y\boxtimes\cal A$-modules on $Y\ts Y$ $$ \cal
B^\prime {\circ}_{\Phi \cal A} \cal B \simeq ( Q {\circ}_{\cal O_X} b(\Phi
\cal A)){\circ}_{\Phi \cal A} \cal B \simeq Q {\circ}_{\cal O_X} \cal B
\simeq Q {\circ}_{\cal O_X} P \circ_{\cal O_{Y}} b(\cal A) \simeq b(\cal A).
$$
Similarly, there is an isomorphism in the derived category of left $\cal
A\boxtimes\cal O_Y$-modules
$$
\cal B^\prime{\circ}_{\Phi\cal A} \cal B \simeq \cal B^\prime{\circ}_{\Phi
\cal A} (b(\Phi \cal A) {\circ}_{\cal O_X} P) \simeq \cal
B^\prime\circ_{\cal O_X}
P \simeq b(\cal A)\circ_{\cal O_{Y}}
Q{\circ}_{\cal O_X} P \simeq b(\cal A).
$$
Moreover, both these isomorphisms
coincide with the following isomorphism of ${\cal O}_{Y\ts Y}$-modules
$$\cal B'\circ_{\Phi\cal
A} \cal B\simeq Q\circ_{\cal O_X} b(\Phi\cal A)\circ_{\Phi\cal A}
b(\Phi\cal A)\circ_{\cal O_X}
P\simeq Q\circ_{\cal O_X} b(\Phi\cal A)\circ_{\cal O_X} P \simeq b(\cal
A).$$ It follows that
$\cal B^\prime{\circ}_{\Phi\cal A} \cal B \simeq b(\cal A)$ in the derived
category of $\cal
A\boxtimes\cal A^{op}$-modules,

Similarly,
$\cal B{\circ}_{\cal A}\cal B'\simeq b(\Phi\cal A)$. It follows that the
compositions $\Phi \Phi': \D^b \cal M(\Phi \cal A) \rightarrow \D^b \cal
M(\Phi \cal A)$ and $\Phi' \Phi: \D^b \cal M(\cal A) \rightarrow \cal D^b
\cal M(\cal A)$ are identity functors.

The composition of $\Phi$ with the forgetting functor can be easily
computed: $$
\cal B {\circ}_{\cal A} M
\simeq ( P {\circ}_{\cal O_Y} b(\cal A)) \circ_{\cal A} M \simeq
P{\circ}_{\cal O_Y} M\ .
$$
Hence, forgetting $\Phi \cal A$-module structure, we just get the transform
with kernel $P$.
\end{pf}

\begin{Rems}
1. In the situation of the theorem if we have another special $D$-algebra
$\cal A^{\prime}$ and a homomorphism of $D$-algebras $\cal A\rightarrow
\cal A^{\prime}$ then we have the corresponding induction and restriction
functors $M\mapsto \cal A^{\prime}\otimes_{\cal A} M$ and $N\mapsto N$
between categories of $\cal A$-modules and $\cal A^{\prime}$-modules. It is
easy to check that the corresponding derived functors commute with our
functors $\Phi$ constructed for $\cal A$ and $\cal A^{\prime}$.

\noindent
2. Let $A$ be an abelian variety, $\hat{A}$ be the dual abelian variety.
Then as was shown in [L] and [R2] the Fourier-Mukai equivalence ${\cal
D}^b(A)\simeq{\cal D}^b(\hat{A})$ extends to an equivalence of the derived
categories of $\cal D$-modules on $A$ and $\O$-modules on the universal
extension of $\hat{A}$ by a vector space. The latter category is equivalent
to the category of modules over the commutative sheaf of algebras $\cal A$
on $\hat{A}$ which is constructed as follows.
Let
$$0\ra \cal O\ra\cal E\ra H^1(\hat{A},\cal O)\otimes\cal O\ra 0$$ be the
universal extension. Then
$\cal A=\Sym(\cal E)/(1_{\cal E}-1)$ where $1_{\cal E}$ is the image of
$1\in\cal O$ in $\cal E$. It is easy to see that $\cal A$ is the dual
special $D$-algebra to the algebra of differential operators on $A$, so our
theorem implies the mentioned equivalence of categories.

\noindent
3. In the case of abelian varieties one can generalize the notion of
special $D$-algebra as follows.
Instead of considering special sheaves on $X\times X$ one can consider
quasi-coherent sheaves on
$X\times X$ admitting filtration with quotients of the form $(\id, t_x)_*L$
where $(\id,
t_x):X\rightarrow X\times X$ is the graph of the translation by some point
$x\in X$, $L$ is a line
bundle algebraically equivalent to zero on $X$. Let us call such sheaves
quasi-special. It is easy
to see that quasi-special sheaves are flat over $X$ with respect to both
projections $p_1$ and
$p_2$, so the operation $\circ$ is exact on them. We can define a
quasi-special algebra as a
quasi-special sheaf $K$ on $X\times X$ together with the associative
multiplication $K\circ
K\rightarrow K$ admitting a unit $\Delta_*\cal O_X\rightarrow K$. Then
there is a Fourier duality
for quasi-special algebras and equivalence of the corresponding derived
categories. The proof of
the above theorem works literally in this situation. Note that modules over
quasi-special algebras
form much broader class of categories than those over special $D$-algebras.
Among these categories
we can find some categories of modules over 1-motives and our Fourier
duality coincides with the
one defined by G. Laumon in [L]. For example, a homomorphism $\phi:\Bbb
Z\rightarrow X$ defines a
quasi-special algebra on $X$ which is a sum of structural sheaves of graphs
of translations by
$\phi(n)$, $n\in\Bbb Z$. The corresponding category of modules is the
category of $\Bbb
Z$-equivariant $\cal O_X$-modules. The Fourier dual algebra corresponds to
the affine group over
$\hat{X}$ which is an extension of $\hat{X}$ by the multiplicative group.
\end{Rems}

\subsection{}
Let $L$ be a Lie algebroid on $Y$ such that $L \simeq \cal O^d_{Y}$ as an $
\cal O_{Y}$ - module. Then for any central extension $\tilde{L}$ of $L$ by
$\cal O_{Y}$, the $D$-algebra \ \ $\cal U^\circ (\tilde{L})$ is special.
Futhermore, one has $\Phi \cal U^\circ (\tilde{L}) \simeq \cal
U^\circ(\tilde{L^\prime})$ for some central extension $ \tilde{L^\prime}$
of a Lie algebroid $L^\prime$ on $X$ such that $L^\prime \simeq \cal O^d_X$
as an $\cal O_X$-module. Indeed, this follows essentially from Lemma
\ref{locallyfree}. One just has to notice that if a $D$-algebra $\cal A$ on
$Y$ has an algebra filtration $\cal A_\bullet$ with $\gr \cal A_\bullet
\simeq S^\bullet (\cal O^d_{Y})$, then $\Phi \cal A$ has an algebra
filtration $\cal F \cal A_\bullet$ with $\gr \Phi \cal A_\bullet \simeq
S^\bullet(\cal O^d_X)$. Note that if $L$ is a successive extension of
trivial bundles then the $D$-algebra $\cal U^\circ(\tilde{L})$ is still
special, but $\Phi \cal U^\circ(\tilde{L})$ is not necessarily of the form
$\cal U^\circ (\tilde{L'})$.

\subsection{}
Assume now that $X$ is an abelian variety. Let $\cal P$ be a Picard
algebroid on $\hat X, \cal D= \cal U^\circ(\cal P)$ be the corresponding
$tdo$. Then $\cal P/ \cal O_{\hat X} \simeq \cal T_{\hat X} \simeq \hat
{\frak g} \otimes_k \cal O _{\hat X}$ is a trivial $\cal O_{\hat
X}$-module, hence $\Phi \cal D \simeq \cal U^{\circ}(\tilde{L^\prime})$ for
some Lie algebroid ${L^\prime}$ on $X$ and its central extension
$\tilde{L^\prime}$ by $\cal O_X$.
\begin{Prop}
Let $\cal D$ be a tdo on $\hat X, \cal P$ be the corresponding Picard
algebroid. Then $\Phi\cal D$ is a tdo on $X$ if and only if the map
$\hat{\frak g} \rightarrow H^1 (\hat X, \cal O)$, induced by the extension
of $\cal O_{\hat X}$-modules $$
0 \rightarrow \cal O_{\hat X}\rightarrow \cal P \rightarrow \hat{\frak g}
\otimes_k \cal O_{\hat X}\rightarrow 0,
$$
is an isomorphism.
\end{Prop}
\begin{pf}
Let $\cal D_\bullet$ be the canonical filtration of $\cal D$. Then $\Phi
\cal D$ is a tdo if only if the class of the extension of $\cal
O_X$-bimodules $$
0 \rightarrow \cal O_X \simeq \Phi \cal D_0 \rightarrow \Phi \cal D_1
\rightarrow
\Phi(\cal D_1/ \cal D_0) \simeq \hat{\frak g} \otimes \cal O_X \rightarrow
0 $$ induces an isomorphism $\hat{\frak g} \otimes \cal O_{\hat X}
\rightarrow \cal T_X$. Thus, it is
sufficient to check that the components of the canonical decomposition $$
\Ext^1_{\cal O_{X \times X}} (\Delta_* \cal O_X, \Delta_*, \cal O_X) \simeq
H^0(X, \cal T)\oplus H^1(X, \cal O_X)
$$
introduced in Lemma \ref{isomorphism}, get interchanged by the
Fourier-Mukai transform, if
we take into
account the natural isomorphisms
$$
H^0(X, \cal T) \simeq \frak g \simeq H^1(\hat X, \cal O), $$ $$
H^1 (X,\cal O_X) \simeq \hat{\frak g} \simeq H^0 (\hat X,\cal T). $$ We
leave this to the reader as a pleasant exercise on Fourier-Mukai transform.
\end{pf}

\subsection{}
Let us describe in more details the data describing a Lie algebroid $L$ on
an abelian
variety $X$ such that $L \simeq V \otimes_k \cal O_X$ as $\cal O_X$-module,
where $V$ is a finite-dimensional $k$-vector space, and a central extension
$\tilde {L}$ of $L$ by $\cal O_X$. First of all, $V=H^0(X,L)$ has a
structure of
Lie algebra, and the structural morphism $L \rightarrow \cal T$ is given by
some $k$-linear map $\beta : V
\rightarrow \frak g = H^0 (X,\cal T)$ which is a homomorphism of Lie
algebras (where $\frak g$ is an
abelian Lie algebra). The central extension $\tilde {L}$ is described (up
to an isomorphism) by a class $\widetilde{\alpha}$ in the first
hypercohomology space space $\Bbb
H^1(X,L^*
\rightarrow \wedge^2 L^* \rightarrow \wedge^3 L^* \rightarrow \dots)$ of
the truncated Koszul complex of $L$. In particular, we have the
corresponding class $\alpha
\in H^1(X, L^*)$, which is just the class of the extension of $\cal
O_X$-modules
$$
0 \rightarrow \cal O_X \rightarrow \tilde{L} \rightarrow L \rightarrow 0\ .
$$ We can consider $\alpha$ as a linear map $V \rightarrow H^1(X, \cal O_X)
= \hat{\frak g}$.
The maps $\alpha$ and $\beta$ get interchanged by the Fourier transform, up
to a sign.

By definition the $D$-algebra associated with $\widetilde{L}$ is a tdo if
and only if $\beta:V\rightarrow\frak g$ is an isomorphism. If an addition
$\alpha:V\rightarrow\hat{\frak g}$ is an isomorphism then the dual
$D$-algebra is also a tdo. Thus, we have a bijection between tdo's with
non-degenerate first Chern class on $X$ and $\hat{X}$ such that the
corresponding derived categories of modules are equivalent. According to
[BB] isomorphism classes of tdo on $X$ are classified by $\Bbb
H^2(X,\Omega^{\ge 1})$ which is an extension of $H^1(X,\Omega^1)\simeq
\Hom(\frak g,\hat{\frak g})$ by
$H^0(X,\Omega^2)=\wedge^2{\frak g}^*$.
Let $U_X\sub \Bbb H^2(X,\Omega^{\ge 1})$ be the subset of elements with
non-degenerate projection to $H^1(X,\Omega^1)$. The duality gives an
isomorphism between $U_X$ and $U_{\hat{X}}$. It is easy to see that under
this isomorphism multiplication by $\lambda\in k^*$ on $U_X$ corresponds to
multiplication by $\lambda^{-1}$ on $U_{\hat{X}}$.

On the other hand, let $\cal A$ be a tdo with trivial $c_1$. In other
words, $\cal A$ corresponds to some global 2-form $\omega$ on $X$. Modules
over $\cal A$ are $\O$-modules equipped with a connection having the scalar
curvature $\omega$. Let $\cal B$ be the dual $D$-algebra on $\hat{X}$ and
let $\wt{L}\ra L=H^0(X,T_X)\otimes\O_{\hat{X}}$ be the corresponding
central extension
of Lie algebroids. We claim that $L$ is just an $\O_{\hat{X}}$-linear
commutative Lie algebra while the central extension $\wt{L}$ is given by
the class
$(e,\omega)\in H^1(L^*)\oplus H^0(\wedge^2 L^*)$, where $e$ is the
canonical element in $H^1(L^*)\simeq H^1(\hat{X},\O)\otimes
H^1(\hat{X},\O)^*$. Indeed, as an $\O_{\hat{X}}$-module $\wt{L}$ is a
universal extension of $H^1(\hat{X},\O)\otimes\O$ by $\O$. Hence, the Lie
bracket defines a morphism of $\O$-modules $\wedge^2 L\ra\wt{L}$. Since
$H^0(\wt{L})=H^0(\O)$ it follows that $[\wt{L},\wt{L}]\subset
\O\subset\wt{L}$. It is easy to see that the Lie bracket is just given by
$\omega:\wedge^2 L\ra\O$.

Recall that the Neron-Severi group of $X$ is identified with
$\Hom^{\sym}(X,\hat{X})\otimes\Bbb Q$ where $\Hom^{\sym}(X,\hat{X})$ is the
group of symmetric homomorphisms $X\ra\hat{X}$. Namely, to a line bundle
$L$ therecorrespondsa symmetrichomomorphism
$\phi_L:X\ra\hat{X}$ sending a point $x$ to $t_x^*L\otimes L^{-1}$ where
$t_x:X\ra X$ is the translation by $x$. One has the natural homomorphism
$c_1:NS(X)\ra \Bbb H^2(X,\Omega^{\ge 1})$ sending a line bundle $L$ to the
class of the ring $D_L$ of differential operators on $L$. For $\mu\in
NS(X)$ we denote by $D_{\mu}$ the corresponding tdo. If $\mu\in NS(X)$ is a
non-degenerate class so that $c_1(\mu)\in U_X$ then the Fourier dual tdo to
$D_{\mu}$ is
\begin{equation}\label{detisom}
\Phi(D_{\mu})=D_{-\mu^{-1}}
\end{equation}
Indeed, it suffices to check this when $\mu$ is a class of a line bundle
$L$, in which case it follows easily from the isomorphism
$$\phi_L^*\det\Phi(L)\simeq L^{-\rk\Phi(L)}$$ and the fact that the dual
tdo to $D_L$ acts on $\Phi(L)$.

Let $E$ be a coherent sheaf which is a module over some tdo on $X$ (then
$E$ is automatically locally free). Following [BB] we say in this case that
there is an integrable projective connection on $E$.

\begin{Prop}
Let $E$ be a vector bundle on $X$ equipped with an integrable projective
connection. Assume that $\det E$ is a
non-degenerate line bundle.
Then $H^i\Phi(E)$ are vector bundles with canonical integrable projective
connections, and the following equality
holds:
$$\phi_{\det E}^*c_1(\Phi(E))
=-\chi(X,E)\cdot\rk E\cdot c_1(E).$$
\end{Prop}
\begin{pf} The first statement follows immediately from the fact that
$\Phi(E)$ is quasi-isomorphic to a complex of modules over the tdo on
$\hat{X}$ dual to $D_{(\det E)^{\frac{1}{r}}}$ where $r=\rk E$. On the
other hand, this tdo acting on $\Phi(E)$ is isomorphic $D_{(\det
\Phi(E))^{\frac{1}{r'}}}$ where $r'=\rk\Phi(E)=\chi(X,E)$. Considering
classes of these dual tdo's and using the isomorphism (\ref{detisom})
applied to $\mu=\frac{1}{r}\phi_{\det E}$ we get the above formula.
\end{pf}

\subsection{} The following two natural questions arise: 1) whether for
every $\mu\in NS(X)$ there exists a vector bundle $E$ on $X$ which is a
module over $D_{\mu}$, 2) what vector bundles on an abelian variety admit
integrable projective connections. To answer these questions we use the
following construction.
Let $\pi:X_1\ra X_2$ be an isogeny of abelian varieties and $E$ be a vector
bundle with an integrable projective connection on $X_1$. Then there is a
canonical integrable projective connection on $\pi_*E$. Indeed, the
simplest way to see this is to use Fourier duality. If $E$ is a module over
some tdo $D_{\lambda}$ on $X_1$ then $\Phi(E)$ is a module over the dual
$D$-algebra $\Phi(D_{\lambda})$ on $\hat{X}_1$. Now we use the formula
$$\pi_*E\simeq \Phi^{-1}\hat{\pi}^*(\Phi(E)),$$ where $\Phi^{-1}$ is the
inverse Fourier transform on $X_2$, hence $\pi_*E$ is a module over
$\Phi^{-1}\hat{\pi}^*\Phi(D_{\lambda})$ which is a tdo on $X_2$.

In particular, the push-forwards of line bundles under isogenies have
canonical integrable projective connections. Also it is clear that if $E$
is a vector
bundle with an integrable projective connection and $F$ is a flat vector
bundle then $E\otimes F$ has a natural integrable projective connection.

Now we can answer the above questions.

\begin{Thm} For every $\mu\in NS(X)$ there exists a vector bundle $E$ which
is a module over $D_{\mu}$.
\end{Thm}
\begin{pf} We can write $\mu=[L]/n$ where $n>0$ is an integer, $[L]$ is a
class of a line bundle $L$ on $X$. Let $[n]_A:A\ra A$ be an endomorphism of
multiplication by $n$. Then $[n]_A^*(\mu)\in NS(X)$ is represented by a
line bundle $L'$. Now we claim that the push-forward $[n]_{A,*}L'$ has a
structure of a module over $D_{\mu}$. Indeed, it suffices to check that
$$c_1([n]_{A,*}L')/\deg([n]_A)=\mu.$$
Let $\Nm_n: NS(X)\ra NS(X)$ be the norm homomorphism corresponding to the
isogeny $[n]_A$. Then the LHS of the above equality is
$\Nm_n([L'])/\deg([n]_A)$. Hence, the pull-back of the LHS by $[n]_A$ is
equal to $[L']=[n]_A^*(\mu)$ which implies our claim. \end{pf}

\begin{Thm} Let $E$ be an indecomosable vector bundle with an integrable
projective connection on an abelian variety $X$. Then there exists an
isogeny of abelian varieties $\pi:X'\ra X$, a line bundle $L$ on $X'$ and a
flat bundle $F$ on $X$ such that $E\simeq \pi_*L\otimes F$. \end{Thm}
\begin{pf}
The main idea is to analyze the sheaf of algebras $A=End(E)$. Namely, $A$
has a flat connection such that the multiplication is covariantly constant.
In other words, it corresponds to a representation of the fundamental group
$\pi_1(X)$ in automorphisms of the matrix algebra. Since all such
automorphisms are inner we get a homomorphism $$\rho:\pi_1(X)\ra PGL(E_0)$$
where $E_0$ is a fiber of $E$ at zero. Now the central extension
$SL(E_0)\ra PGL(E_0)$ induces a central extension of $\pi_1(X)=\Z^{2g}$ by
the group of roots of unity of order $rk E$. This central extension splits
on some subgroup of finite index $H\sub\pi_1(X)$. In other words, the
restriction of $\rho$ to $H$ lifts to a homomorphism $\rho_H:H\ra GL(E_0)$.
Let $\pi:\wt{X}\ra X$ be an isogeny corresponding to $H$, so that $\wt{X}$
is an abelian variety with $\pi_1(\wt{X})=H$. Then $\rho_H$ defines a flat
bundle $\wt{F}$ on $\wt{X}$ such that $$\pi^*A\simeq End(\wt{F})$$ as
algebras with connections. It follows that $\pi^*E\simeq L\otimes\wt{F}$
for some line bundle $L$ on $\wt{X}$. Thus, $L$ is a direct summand of
$\pi_*(L\otimes\wt{F})$.
Note that there exists a flat bundle $F$ on $X$ such that $\wt{F}\simeq
\pi^*F$ (again the simplest way to see this is to use the Fourier duality).
Hence, $L$ is a direct summand of $\pi_*L\otimes F$. It remains to check
that all indecomposable summands of the latter bundle have the same form.
This follows from the following lemma.
\end{pf}

\begin{Lem} Let $\pi:X_1\ra X_2$ be an isogeny of abelian varieties, $L$ be
a line bundle on $X_1$, $F$ be an indecomposable flat bundle on $X_1$.
Assume that $\pi_*(L\otimes F)$ is decomposable. Then there exists a
non-trivial factorization of $\pi$ into a composition
$$X_1\stackrel{\pi'}{\ra}X'_1\ra X_2$$
such that
$L\simeq (\pi')^*L'$ for some line bundle $L'$ on $X'_1$. \end{Lem}
\begin{pf} By adjunction and projection formula we have
$$\End(\pi_*(L\otimes F))\simeq \Hom(\pi^*\pi_*(L\otimes F), L\otimes
F)\simeq\oplus_{x\in K}\Hom(t_x^*L\otimes F,L\otimes F)$$ where $K\subset
X_1$ is the kernel of $\pi$. If $t_x^*L\simeq L$ for some $x\in K$, $x\neq
0$ then $L$ descends to a line bundle on the quotient of $X_1$ by the
subgroup generated by $x$. Otherwise, we get $\End(\pi_*(L\otimes
F))\simeq\End(F)$, hence, $\pi_*(L\otimes F)$ is indecomposable. \end{pf}

\section{Noncommutative \'Etale morphisms} \label{et-sec}

This section provides the setting we will need to discuss microlocalization.

\subsection{}
\begin{Def} (cf. [K]) A ring homomorphism $A'\morph{\phi}A$ is called a central
extension if $\phi$ is surjective, $ker(\phi)$ is a central ideal and
$ker(\phi)^2=0$.
\end{Def}

\begin{Def}\label{formally etale} Let $\rings$ denote the category of
asssociative rings.
Let
$\cal C\subset\rings$ be a full subcategory. A morphism $R\morph{\a}S$ in
$\cal C$
is formally \'etale if, for every
commutative diagram $\alpha, \beta, \gamma, \delta$ \vskip 2pt
\begin{equation}\label{etale morphism}
\begin{array}{ccc} R & \lrar{\delta} & A' \\ \ldar{\alpha}
&\lbrurar{\epsilon} & \ldar{\gamma} \\ S & \lrar{\beta} & A \end{array}
\end{equation}
in $\cal C$, with $\gamma$ a central extension, there exists a unique morphism
\dsp{S\morph{\epsilon}A'} such that diagram (\ref{etale morphism})
commutes. \end{Def}

\begin{Ex}\label{standard ring example} Let $R$ be a ring and let
$a_0,a_1,....,a_n\in R$.
Let
$S$ be the
$R$-algebra generated by elements $z,u$ subject to the relations
\begin{equation}\sum a_iz^i=0\ ,\ u\sum i a_iz^{i-1}=1\ ,\ \sum i
a_iz^{i-1}u=1\
.\end{equation} Then the natural map \dsp{R\morph{\a}S} is formally \'etale
in $\rings$.\end{Ex}
\begin{pf} Consider a commutative diagram \vskip 2pt
\begin{equation}
\begin{array}{ccc} R & \lrar{\delta} & A' \\ \ldar{\alpha} && \ldar{\gamma}
\\ S & \lrar{\beta} & A \end{array}
\end{equation}
as in definition \ref{formally etale}. Let $I=ker(\gamma)$. Choose
$x,y\in A'$ such that $\gamma(x)=\b(z)$, $\gamma(y)=\b(u)$. It must be
shown that there
are unique elements $p,q\in I$ such that \begin{equation}
\begin{matrix}
\sum \delta(a_i)(x+p)^i&=0\ ,\\
\ (y+q)\sum i \delta(a_i)(x+p)^{i-1}&=1\ ,\\ \ \sum i
\delta(a_i)(x+p)^{i-1}(y+q)&=1\ .\end{matrix}\end{equation} Since $I^2=0$,
the equations are uniquely solved by setting \begin{equation}
p=-y\sum \delta(a_i)x^i\ \ ,\ \ q=y(1-\sum i\delta(a_i)(x+p)^{i-1} y)\
.\end{equation}
\end{pf}

\subsection{}
Let us recall some definitions from [K]. For any associative algebra $R$, the
NC-filtration on
$R$ is the decreasing filtration $\{F^d R\}_{d\ge 0}$ defined by setting
$$F^d R=\sum_{i_1+\ldots+i_m=d} R\cdot R_{i_1} \cdot R \cdot \ldots \cdot
R\cdot R_{i_m}
\cdot R$$ where $R_0=R$, $R_{i+1}=[R,R_i]$ are the terms of the lower
central series for $R$
considered as a Lie algebra (we use a different indexing from Kapranov's).
This filtration
is compatible with multiplication and the associated graded algebra is
commutative.

The category $\NN_d$ is the category of associative algebras $R$ with
$F^{d+1}R=0$. For
example,
$\NN_0$ is the category of commutative algebras. For every $d$ there is a
pair of adjoint
functors
$r_d:\NN_d\ra\NN_{d-1}$ and $i_d:\NN_{d-1}\ra\NN_d$, where $i_d$ is the
natural inclusion,
$r_d(R)=R/F_dR$. Note that if $R\in\NN_d$ then $F^dR\sub R$ is a central
ideal with zero
square. Thus, $R$ is a central extension of $r_d(R)$. Indeed, $\NN_d$ is
the category of rings $A$ which are obtained as the composition of $d$
central extensions,
$$
A\ra A_1\ra A_2\ra ...\ra A_d$$
with $A_d$ commutative.

\begin{Lem}\label{der} Let $R\ra S$ be a formally \'etale morphism in
$\NN_d$, $M$ be an
$S^{ab}$-module. Then the natural map $\Der(S,M)\ra\Der(R,M)$ is a
bijection. \end{Lem}

\begin{pf} Given a central $S$-bimodule $M$ we can define a trivial central
extension of $S$
by
$M$: $\wt{S}=S\oplus M$. Then derivations from $S$ to $M$ are in bijective
correspondence
with splittings of the projection $\wt{S}\ra S$. Hence, the assertion. \end{pf}

\begin{Prop}\label{sasha's nice observation} Let \dsp{R\morph{\alpha}S} be
a formally \'etale
morphism in
$\NN_d$. Let
\begin{equation}
\begin{array}{ccc} R & \lrar{\delta} & A' \\ \ldar{\alpha} & &
\ldar{\gamma} \\ S & \lrar{\beta} & A \end{array}
\end{equation} be a commutative diagram in $\rings$, such that $\beta$ is
surjective and $\gamma$ is a central extension. Then $A'\in\NN_d$.\end{Prop}

\begin{pf}
Note that apriori we know from this diagram that $A\in\NN_{d}$, hence,
$A'\in\NN_{d+1}$. We want to prove that $F^{d+1}A'=0$. Since $F^{d+2}A'=0$
it suffices to
prove that for every sequence of positive numbers $i_1,\ldots,i_m$ such
that $i_1+\ldots+i_m=d+1$ one has
$$A'_{i_1}\cdot A'_{i_2}\cdot
\ldots\cdot A'_{i_m}=0.$$ We use descending induction in $m$. Assume that
this is true for
$m+1$. Then we can define a map
\begin{align*} &D:S^{m+d+1}\ra I: \\ &(s_1,\ldots,s_{m+d+1})\mapsto
[a'_1,[\ldots,[a'_{i_1},a'_{i_1+1}]\ldots]]\cdot
[a'_{i_1+2},[\ldots,[a'_{i_1+i_2+1},a'_{i_1+i_2+2}]\ldots]]]\cdot \ldots
\end{align*} where $a'_i\in A'$ are such that $\gamma(a'_i)=\beta(s_i)$.
This map is
well-defined since
$a'_i$ are well-defined modulo $I$ which is a central ideal. Now the
induction assumption
implies that $D$ is a derivation in every argument. Hence, applying Lemma
\ref{der} we
conclude that $D=0$.
\end{pf}

\begin{Thm}\label{enlarging the category} Let \dsp{R\morph{\alpha}S} be a
formally \'etale
morphism in
$\NN_d$. Then \dsp{R\morph{\alpha}S}is a formally \'etale morphism in
$\rings$.\end{Thm}
\begin{pf} This follows easily from propostion \ref{sasha's nice
observation}.\end{pf}

Let $\NN_{\infty}$ denote the category of rings that are complete with
respect to the
NC-filtration.

\begin{Thm} Let \dsp{R\morph{\alpha}S} be formally \'etale in
$\NN_{\infty}$, with
$R\in\NN_d$. Then $S\in\NN_d$.\end{Thm}

\begin{pf} The natural morphism \dsp{R\ra r_{d+i}(S)} is formally \'etale
in $\NN_{d+i}$ for all $i\ge 0$. Proposition \ref{sasha's nice observation}
applied to the
diagram
\begin{equation}
\begin{array}{ccc} R & \lrar{} & r_{d+i+1}(S) \\ \ldar{} & & \ldar{} \\
r_{d+i}(S) & \lrar{=} & r_{d+i}(S) \end{array}
\end{equation}
shows that $F^{d+i}(S)\subset F^{d+i+1}(S)$ for all $i\ge 1$. Hence the
assertion.\end{pf}

\begin{Ex}\label{standard example in Nd} Let $R$ be a ring belonging
$\NN_d$ for some $d$.
Let
\dsp{R\morph{\a}S} be as in example \ref{standard ring example}. Let $\hat
S$ denote
the completion of $S$ with respect to the NC-filtration. Then $\hat S$
belongs to $\NN_d$
and the
natural morphism $R\ra \hat{S}$ is formally \'etale (in $\rings$). As in
the commutative case, we call such a morphism {\em standard}.\end{Ex}

\subsection{}
The category $\NCS^d$ of NC-schemes of degree $d$ (in Kapranov's
terminology ``NC-nilpotent of
degree $d$") is constructed in the same way as the commutative category of
schemes using
$\NN_d$ instead of $\NN_0$ as coordinate rings of affine schemes. For a
scheme $X\in \NCS^d$
we denote by
$\NCS^d_{/X}$ the category of NC-schemes of degree $d$ over $X$. The
morphisms in
$\NCS^d_{\/X}$ are denoted by $\Hom_X(\cdot,\cdot)$.

As in affine case we have natural adjoint functors $r_d:\NCS^d\ra
\NCS^{d-1}$ and
$i_d:\NCS^{d-1}\ra \NCS^d$. In particular, we have the abelianization
functor $\NCS^d\ra
\NCS^0:X\mapsto X^{ab}$ given by the composition $r_1r_2\ldots r_d$.

A morphism $Z\ra\wt{Z}$ of NC-schemes of degree $d$ is called a nilpotent
thickening if it
induces an isomorphism of underlying topological spaces and
$\O_{\wt{Z}}\ra\O_Z$ is a
surjection with nilpotent kernel.

\begin{Def} A morphism $Y\ra X$ in $\NCS^d$ is called formally smooth
(resp. formally
unramified) if for every nilpotent thickening $Z\sub \wt{Z}$ in
$\NCS^d_{/X}$ the map
$\Hom_X(\wt{Z},Y)\ra\Hom_X(Z,Y)$ is surjective (resp. injective). A
morphism is called
formally
\'etale if it is both formally smooth and formally unramified. A morphism
is called \'etale
if it is \'etale and the corresponding morphism of commutative schemes
$Y^{ab}\ra X^{ab}$ is
locally of finite type. \end{Def}

\begin{Prop}\label{gen} Let $P$ be a property of being formally smooth
(resp. formally
unramified, resp. formally \'etale). \newline a) Let $f:Y\ra X$ be a
morphism in $\NCS^d$
with property $P$. Then the same property holds for $r_d(f)$. \newline b)
If $f:Y\ra X$ and
$g:Z\ra Y$ are morphisms in $\NCS^d$ having property $P$ then $f\circ g$
also has this
property. \newline c) If $f:Y\ra X$ is formally unramified morphism in
$\NCS^d$, $g:Z\ra Y$
is a morphism in $\NCS^d$ such that $f\circ g$ has property $P$ then $g$
also has this
property. \newline d) An open morphism $U\ra X$ is \'etale. \newline e) A
morphism $f:Y\ra X$ is \'etale if and only if there exists an open
covering
$X=\cup X_i$ and for every $i$ an open covering $Y_{ij}$ of $f^{-1}(X_i)$
such that all the
induced morphisms
$Y_{ij}\ra X_i$ are \'etale.
\end{Prop}

The proof is straightforward.

Theorem \ref{enlarging the category} has the following global version.

\begin{Thm} Let $f$ be a formally \'etale morphism in $\NCS^{d-1}$. Then
$i_d(f)$ is a
formally
\'etale morphism in $\NCS^d$. \end{Thm}

Now we observe that the topological invariance of \'etale morphisms remains
valid in the
present context.

\begin{Thm} For any $X\in \NCS^d$ the canonical functor $Y\mapsto Y^{ab}$
from the category
of
\'etale $X$-schemes to that of \'etale $X^{ab}$-schemes is an equivalence.
\end{Thm}

\begin{pf} First we claim that the functor in question is fully faithful.
Indeed, let $Y_1,
Y_2\in\NCS^d$ be \'etale $X$-schemes. Then since $Y_1^{ab}\ra Y_1$ is a
nilpotent thickening
and
$Y_2\ra X$ is \'etale the natural map
$$\Hom_X(Y_1,Y_2)\ra\Hom_X(Y_1^{ab},Y_2)\simeq
\Hom_{X^{ab}}(Y_1^{ab},Y_2^{ab})$$ is an isomorphism as required.

To prove the surjectivity of the functor it suffices to do it locally. Thus
we may assume
that the morphism $Y^{ab}\ra X^{ab}$ is a standard \'etale extension of
commutative rings
$R^{ab}\ra S^{ab}$ where $S^{ab}=(R^{ab}[z_0]/(f_0(z_0)))_{f'_0(z_0)}$,
$f_0$ is a unital polynomial. But such a morphism lifts to a standard
\'etale morphism in
$\NN_d$, as in example \ref{standard example in Nd}. \end{pf}

\begin{Cor}
A morphism $f:Y\ra X$ in $\NCS^d$ is \'etale if and only if there exists an
open
covering
$Y=\cup Y_i$ such that all the
induced morphisms
$Y_{i}\ra f(Y_i)$ are standard \'etale morphisms. \end{Cor}

\section{Microlocalization}

\subsection{}
Let us return to the setting of Theorem (\ref{main1}). We assume that the
$D$-algebra $\cal A$ is equipped with an increasing algebra filtration $\cal
A_{\bullet}$ (so $\cal A_i \cal A_j\subset \cal A_{i+j}$) such that $\cal
A_{-1}=0$, $\cal
A_0=\cal O_Y$, the associated graded algebra $\gr(\cal A)$ is commutative
and is generated by
$\gr(\cal A)_1$ over $\cal O_Y$. In particular, the left and right actions
of $\cal O_Y$ on
$\gr(\cal A)_i$ are the same.
We will call such a filtration {\it special} if there exists a sheaf of
flat, commutative graded
$\cal O_S$-algebras $C$, generated over $\cal O_S$ by $C_1$, and an
isomorphism of
graded algebras
$\gr(\cal A)\simeq {\pi_S^Y}^*( C)$. By lemma \ref{sheaves on S}, such an
isomorphism induces
an isomorphism
$\gr(\Phi\cal A)\simeq {\pi_S^X}^*( C)$.

Given a $D$-algebra $\cal A$ with a special filtration $\cal A_{\bullet}$
we can form the
corresponding sequence of graded algebras $\gr_{(n)}(\cal A)$ for $n\ge 0$
by setting $$\gr_{(n)}(\cal A)=\oplus_{i=0}^{\infty} \cal A_i/\cal
A_{i-n-1}.$$ In particular,
$\gr_{(0)}(\cal A)=\gr(\cal A)$ is commutative while for $n\ge 1$ there is
a central element $t$ in
$\gr_{(n)}(\cal A)_1=\cal A_1$ (corresponding to $1\in \cal A_0\subset \cal
A_1$) such that $t^{n+1}=0$ and $\gr_{(n)}(\cal A)/(t)=\gr(\cal A)$.
These algebras form a projective system via the natural projections
$\gr_{(n+1)}(\cal A)\rightarrow \gr_{(n)}(\cal A)$.

Consider for each $n$ the NC-scheme
$\PP_{n}(\cal A)=\Proj(\gr_{(n)}(\cal A))$ corresponding to $\gr_{(n)}(\cal
A)$ via the noncommutative analogue of $\Proj$ construction. We denote by
$\cal
D^-(\PP_{n}(\cal A))$ the (bounded from above) derived category of left
quasi-coherent sheaves
$\PP_{n}(\cal A)$. Similar to the commutative case there is a natural
localization functor
$M\mapsto\widetilde{M}$ from the category of graded $\gr_{(n)}(\cal
A)$-modules to the category of
quasi-coherent sheaves on $\PP_{n}(\cal A)$.

If $M$ is a left $\cal A$-module (quasi-coherent over ${\cal O}_Y$)
equipped with an increasing module filtration $M_{\bullet}$ then for every
$n>0$ we can form the
corresponding graded $\gr_{(n)}(\cal A)$-module $\oplus M_i/M_{i-n}$, hence
the corresponding
quasi-coherent sheaf $$\ml_{n}(M)=(\oplus_i M_i/M_{i-n})^{\widetilde{}}.$$

The above NC-schemes are connected by a sequence of closed embeddings
$i_n:\PP_n(\cal A)\hookrightarrow \PP_{n+1}(\cal A)$ and the quasi-coherent
sheaves
$\ml_n(M)$ satisfy $\ml_n(M)=\ml_{n+1}(M)|_{\PP_n(\cal A)}$. In other
words, the system
$(\ml_n(M))$ corresponds to a quasi-coherent sheaf on the formal NC-scheme
$\PP_{\infty}(\cal
A)=\operatorname{inj.}\lim\PP_n(\cal A)$.

Our aim now is to establish an equivalence of derived categories of sheaves
on $\PP_n(\cal A)$ and $\PP_n(\Phi \cal A)$. We will prove something
stronger -- namely
that such an equivalence exists \'etale locally on $Proj(C)$.

With begin a Zariski local version. Clearly, $\circ$ commutes with flat
base change on $S$, so
we may assume $S$ is affine.
The isomorphisms
$\gr(\cal A)\simeq {\pi_S^Y}^*(\cal C)$ and $\gr(\Phi\cal A)\simeq
{\pi_S^X}^*(\cal C)$ give us isomorphisms $\PP_1(\cal A)=Y\ts\Proj(C)$
and $\PP_1(\Phi\cal A)= X\ts\Proj(C)$.
Let $f$ be a section of $C_1$.
It defines a Zariski open subset $D_f\subset\Proj(C)$ which is the spectrum of
$C_{(f)}$, the degree zero part in the localization of $C$ by $f$. Hence,
we have the
corresponding open subset
$Y\ts D_f\sub\PP_0(\cal A)$. Since $\PP_n$ has the same underlying
topological space as $\PP_0$
we have the corresponding open subscheme $\PP_n(\cal A)_f\subset\PP_n(\cal
A)$ for every $n\ge 0$.
We claim that we can identify $\PP_n(\cal A)_f$ with the spectrum of some
sheaf of algebras over
$Y$. Namely, we have the surjection
$\cal A_1\ra\gr(\cal A)_1$ and locally we can lift $f$ to a section
$\wt{f}\in\cal A_1$.
Consider the graded $\O_Y$-algebra
$$\gr_{(n)}(\cal A)_{f}:=\gr_{(n)}(\cal A)_{\wt{f}}.$$ It is easy to see
that this algebra
doesn't depend on a choice of the lifting element $\wt{f}$ (since two
liftings differ by a
nilpotent), hence this algebra is defined globally over $Y$. Now let
$\gr_{(n)}(\cal A)_{(f)}$ be
the degree zero component in $\gr_{(n)}(\cal A)_f$. Then
$\Spec(\gr_{(n)}(\cal A)_{(f)})$ is an
open subscheme in $\PP_n(\cal A)$ with the underlying open subset $Y\ts
D_f$, hence
$\Spec(\gr_{(n)}(\cal A)_{(f)})=\PP_n(\cal A)_f$.

For a graded $\gr_{(n)}(\cal A)$-module $M$ we have the corresponding
quasi-coherent
sheaf $\wt{M}$ on $\PP_n(\cal A)$. The restriction of $\wt{M}$ to
$\PP_n(\cal A)_f$ is the sheaf
associated with $\gr_{(n)}(\cal A)_{(f)}$-module $M_{(f)}$, the degree zero
part in the
localization of $M$ with respect to some local lifting of $f$.

Let as called a graded sheaf {\it graded special} if every of its graded
components is a special sheaf.

\begin{Lem}\label{lemloc} For every element $f\in C_1$ the algebra
$\gr_{(n)}(\cal A)_f$ is a graded special $D$-algebra on $Y$. There is a
canonical
isomorphism of graded algebras \begin{equation}\label{Philoc}
\Phi(\gr_{(n)}(\cal A)_f)\simeq\gr_{(n)}(\Phi\cal A)_f. \end{equation}
\end{Lem}
\begin{pf}
First of all the Lemma is obvious for $n=0$: in this case $\gr(\cal
A)_f\simeq {\pi_S^Y}^*
(C_f)$ and
$$\Phi(\gr(\cal A)_f)\simeq \gr(\Phi\cal A)_f\simeq {\pi_S^X}^* (C_f).$$
Now for
$n>0$ consider the filtration on $\gr_{(n)}(\cal A)$ by two-sided principal
ideals $I^k=(t^k)$
where $t\in\gr_{(n)}(\cal A)_1$ is the central element corresponding to
$1\in\cal A_1$. Then
$I^0=\cal A$, $I^n=0$, and $I^k/I^{k+1}\simeq\gr(\cal A)(-k)$ as
$\gr_{(n)}(\cal A)$-module for
$0\le k<n$. Localizing this filtration we get a filtration by two
sided-ideals $I^k_f$ in
$\gr_{(n)}(\cal A)_f$ with associated graded quotients $\gr(\cal A)_f(-k)$.
Thus, $\gr_{(n)}(\cal
A)_f$ is graded special.

To construct an isomorphism (\ref{Philoc}) we notice that since
$\Phi(\gr_{(n)}(\cal A)_f)$ is a nilpotent extension of $\Phi(\gr(\cal
A)_f)\simeq
\gr(\Phi\cal A)_f$, any local lifting of $f$ is invertible in
$\Phi(\gr_{(n)}(\cal A)_f)$.
Therefore, by universal property we get a homomorphism $$\gr_{(n)}(\Phi\cal
A)_f\ra
\Phi(\gr_{(n)}(\cal A)_f).$$
Using the above filtration on $\gr_{(n)}$ one immediately checks that
this is an isomorphism.
\end{pf}

\begin{Thm}\label{main2}
Assume that $H^0(X,\cal O_X)=H^0(Y,\O_Y)=\C$. Then for every $n\ge 0$
there is a canonical equivalence of categories $$\Phi_{(n)}:\cal
D^-(\PP_{n}(\cal A))\rightarrow \cal D^-(\PP_{n}(\Phi \cal A))$$
commuting with functors $i_{n*}$ and $i_n^*$. Moreover, assume that $M$ is
a left $\cal A$-module with an increasing module filtration such that
for some integer $d$ and for $i$ sufficiently large $\Phi(M_i/M_{i-1})$ is
concentrated in degree
$d$ (as an object of $\cal D(X)$). Then
$$\ml_n(\Phi(M))\simeq\Phi_{(n)}(\ml_n(M)).$$ \end{Thm} \begin{pf}
The proof of this theorem is similar to the proof of Theorem \ref{main1}.
First we note that the definition of the operation $\circ$ from
\ref{circle} works for
non-commutative schemes as well (the only difference is that now whenever
we need to take the
opposite $D$-algebra we have to pass to the opposite scheme as well). Now
we just want to
construct some quasi-coherent sheaves (perhaps shifted) on $\PP_n(\cal
A)\times\PP_n(\Phi\cal
A^{op})$ and on $\PP_n(\Phi\cal A)\times\PP_n(\cal A^{op})$ such that both
their
$\circ$-composition are equal to the
structure sheaves of diagonals (note that although there is no embedding of
a noncommutative scheme $Z$ into $Z\times Z^{op}$ we still can define an
analogue
of the structure sheaf of the diagonal $\delta_Z$ which is a quasi-coherent
sheaf on $Z\times
Z^{op}$). The definition of these sheaves is the following. First we
observe that
$\PP_n(\Phi\cal A)\times \PP_n(\cal A^{op})=\Proj(\cal A_{XY})$ where $\cal
A_{XY}$ is the following graded algebra on $X\times Y$: $$\cal
A_{XY}=\oplus_i \gr_{(n)}(\Phi\cal A)_i\boxtimes \gr_{(n)}(\cal A^{op})_i.$$
Next we remark that the sheaf
$$\cal B=P\circ_{\cal O_Y} b(\cal A)=b(\Phi\cal A)\circ_{\cal O_X} P$$
introduced
in the proof of Theorem \ref{main1} has a natural filtration $$\cal
B_i=P\circ_{\cal O_Y} b(\cal A)_i= b(\Phi\cal A)_i\circ_{\cal O_X} P.$$ It
follows that the sheaf
$\oplus_i\cal B_{2i}/\cal B_{2i-n}$ has a natural structure of graded $\cal
A_{XY}$-module, so we can set $$\cal B_{ml}=(\oplus_i \cal B_{2i}/\cal
B_{2i-n})^{\widetilde{}}$$ which is a quasi-coherent sheaf on
$\PP_n(\Phi\cal A)\times \PP_n(\cal
A^{op})$. Similarly, one can define the quasi-coherent sheaf $\cal B'_{ml}$
on $\PP_n(\cal
A)\times\PP_n(\Phi\cal A^{op})$ starting with the sheaf $\cal
B'=Q\circ_{\cal O_X} b(\Phi\cal A)$.

It remains to compute
$\cal B_{ml}\circ_{\cal O_{\PP_n(\cal A)}}\cal B'_{ml}$ and $\cal
B'_{ml}\circ_{\cal O_{\PP_n(\Phi\cal A)}}\cal B_{ml}$. The idea is the
following: we cover
$\Proj(C)$ by affine subsets $D_f$, where $f$ runs through $C_1$. For every
$f\in C_1$ we'll construct a canonical isomorphism between restrictions of
$\cal B_{ml}\circ_{\cal O_{\PP_n(\cal A)}}\cal B'_{ml}$ and the structure
sheaf of diagonal in
$\PP_n(\Phi\cal A)\times \PP_n(\Phi\cal A^{op})$ to the open subscheme
$\PP_n(\Phi\cal A)_f\times\PP_n(\Phi\cal A^{op})_f$. These isomorphisms
will be compatible on intersections, so they will glue into a global
isomorphism.

The following notation will be useful: for a sheaf $\cal F$ on one of our
schemes and an element $f\in C_1$ we denote by $\cal F_f$ the restriction
of $\cal F$ to the open subscheme defined by $f$.

Under identification of the underlying topological space of $\PP_n(\Phi\cal
A)\times\PP_n(\cal A^{op})$ with $X\times\Proj(C)\times Y\times\Proj(C)$
the support of $\cal B_{ml}$ is the diagonal $X\times Y\times\Proj(C)$.
Using this fact it is fairly easy to see that $$(\cal B_{ml}\circ_{\cal
O_{\PP_n(\cal A)}}\cal B'_{ml})_f= \cal B_{ml,f}\circ_{\cal O_{\PP_n(\cal
A)_f}}\cal B'_{ml,f}$$ It remains to compute the $\circ$-composition in the
RHS. This is easier than the original problem because the sheaf $\cal
B_{ml,f}$ (resp. $\cal B'_{ml,f}$) live on affine schemes over $X\times Y$
(resp. $Y\times X$). Namely, $$\PP_n(\Phi\cal A)_f\times\PP_n(\cal
A^{op})_f= \Spec(\gr_{(n)}(\Phi\cal A)_{(f)}\boxtimes\gr_{(n)}(\cal
A^{op})_{(f)}).$$

According to Lemma \ref{lemloc} we have dual $D$-algebras $\gr_{(n)}(\cal
A)_{(f)}$ on $Y$
and $\gr_{(n)}(\Phi\cal A)_{(f)}$ on $X$. Hence,
we can apply Theorem \ref{main1} to these $D$-algebras. Let us denote by
$\cal B_{(f)}$ the
$\gr_{(n)}(\Phi\cal A)_{(f)}\boxtimes\gr_{(n)}(\cal A^{op})_{(f)}$-module
constructed in the proof
of the cited theorem (where it is called $\cal B$). We claim that there is
a canonical isomorphism
of the $\cal B_{ml,f}$ with the sheaf on $\Spec(\gr_{(n)}(\Phi\cal
A)_{(f)}\boxtimes\gr_{(n)}(\cal
A^{op})_{(f)})$ obtained by localization of $\cal B_{(f)}$. This claim
(together with an easy check
of the compatibility of isomorphisms on intersections) would allow to
finish the proof by referring
to Theorem \ref{main1}. It remains to construct an isomorphism between the two
$\gr_{(n)}(\Phi\cal A)_{(f)}\boxtimes\gr_{(n)}(\cal A^{op})_{(f)}$-modules:
$\cal B_{(f)}$ and
$(\oplus_i\cal B_{2i}/\cal B_{2i-n})_{(f\otimes f)}$ (the localization of
the latter module is
clearly $\cal B_{ml,f}$). Recall that by definition
$\cal B_{(f)}=P\circ_{\cal O_Y}\gr_{(n)}\cal A_{(f)}$. Also, it is clear
that $(\oplus_i\cal B_{2i}/\cal B_{2i-n})_{(f\otimes f)}$ is isomorphic to
the degree zero part in the localization of $$\gr_{(n)}(\cal
B)=\oplus_i\cal B_i/\cal B_{i-n}= P\circ_{\cal O_Y}\gr_{(n)}(\cal A)$$ by
$\wt{f}\otimes 1$ and $1\otimes \wt{f'}$, where $\wt{f}$ is a local lifting
of $f$ to $\Phi\cal A_1$, $\wt{f'}$ is a local lifting of $f$ to $\cal
A_1$. Thus, it suffices to construct a graded isomorphism between
$P\circ_{\cal O_Y}\gr_{(n)}(\cal A)_f$ and $\gr_{(n)}(\cal
B)_{\wt{f}\otimes 1,1\otimes \wt{f'}}$. According to Lemma \ref{lemloc} we
have
$$P\circ_{\cal O_Y}\gr_{(n)}(\cal A)_f\simeq \gr_{(n)}(\Phi\cal
A)_f\circ_{\cal O_X} P$$ so the assertion follows. \end{pf}

Note that we have
canonical invertible $\O_{\PP_n}$-bimodules on $\PP_n(\cal A)$:
$$\O_{\PP_n}(m)=\ml_n(\gr_{(n)}(\cal A)(m))$$ where $M\mapsto M(m)$ denotes
the shift of grading. In particular, we have the automorphism $$M\mapsto
M(1)=\O(1)\otimes_{\O} M$$
of the category $\D^-(\PP_n(\cal A))$. It is easy to see that the above
equivalence respects these automorphisms.

\subsection{} One can generalize Theorem \ref{main1} to the case of
NC-schemes of finite degree. Namely, there is a natural notion of support
of a quasi-coherent sheaf on such a scheme (just the support of the
corresponding sheaf on the reduced commutative scheme), hence, the
definition of $D$-algebra makes sense. Now the proof of Theorem \ref{main1}
works almost literally in this case. Moreover, it seems plausible for the
NC-schemes $\PP_n(\cal A)$ one can consider slightly more general
$D$-algebras than special ones. Namely, instead of requiring the existence
of filtration with graded factors isomorphic to $\O$ it suffices to require
the existence of filtration with factors $\O^{ab}$, plus one should require
$D$-algebra to be flat as left and right $\O$-module.

\subsection{\'Etale local version of the equivalence}

Let $U\ra Z$ be an \'etale morphism of $S$-schemes. Then we have the
corresponding \'etale morphism $Y\times U\ra\PP_1(\cal A)$. By topological
invariance of \'etale
category for every $n\ge 1$ this morphism extends to an \'etale morphism of
NC-schemes
$$j:\PP_n(\cal A)_U\ra\PP_n(\cal A).$$
Similarly we have an NC-scheme $\PP_n(\Phi\cal A)_U$, and an \'etale
morphism $j:\PP_n(\Phi\cal A)_U\ra\PP_n(\cal A)$.

\begin{Thm} In the above situation the categories $\cal D^-(\PP_n(\cal
A)_U)$ and
$\cal D^-(\PP_n(\Phi\cal A)_U)$ are canonically equivalent. \end{Thm}
\begin{pf} Recall that in the proof of Theorem \ref{main2} we have
constructed a quasi-coherent sheaf (up to shift) $\cal B_{ml}$ on
$\PP_n(\Phi\cal A)\times \PP_n(\cal A^{op})$ supported on the diagonal
$\Delta_Z:X\times Y\times Z\hra X\times Z\times Y\times Z$. Moreover, the
restriction of $B_{ml}$ to $\PP_1(\Phi\cal A)\times\PP_n(\cal A^{op})$ is
actually obtained from the sheaf $P$ on $X\times Y$ via first pulling back
to $X\times Y\times Z$ and then pushing forward by $\Delta_Z$. We have the
following diagram of \'etale morphisms of NC-schemes: \begin{equation}
\begin{array}{ccc}
\PP_n(\Phi\cal A)_U\times \PP_n(\cal A^{op})_U &\lrar{\id\times j}&
\PP_n(\Phi\cal A)_U\times \PP_n(\cal A^{op})\\ \ldar{j\times\id} & &
\ldar{j\times\id}\\ \PP_n(\Phi\cal A)\times \PP_n(\cal A^{op})_U
&\lrar{\id\times j} &\PP_n(\Phi\cal A)\times \PP_n(\cal A^{op}) \end{array}
\end{equation}
Now we claim that there exists a quasi-coherent sheaf $\cal B_{ml,U}$ on
$\PP_n(\Phi\cal A)_U\times \PP_n(\cal A^{op})_U$ supported on the diagonal
$\Delta_U:X\times Y\times U\hra X\times U\times Y\times U$ such that
\begin{equation}\label{iso1}
(\id\times j)_*\cal B_{ml,U}\simeq (j\times\id)^*\cal B_{ml} \end{equation}
\begin{equation}\label{iso2}
(j\times\id)_*\cal B_{ml,U}\simeq (\id\times j)^*\cal B_{ml}.
\end{equation} and such that the restriction of $\cal B_{ml,U}$ to
$\PP_1(\Phi\cal A)_U\times\PP_1(\cal A)$ is isomorphic to
$\Delta_{U,*}(p_{XY}^*P)$. Indeed, consider the quasi-coherent sheaf
$(j\times j)^*\cal B_{ml}$ on $\PP_n(\Phi\cal A)_U\times \PP_n(\cal
A^{op})_U$. It is is supported on $(j\times j)^{-1}(X\times Y\times Z)$
where $X\times Y\times Z$ is the relative diagonal in $X\times Y\times
Z\times Z$. Now since $j$ is \'etale, the relative diagonal $X\times
Y\times U$ is a connected component in $(j\times j)^{-1}(X\times Y \times
Z)$. Now we just set $\cal B_{ml,U}$ to be the direct summand of $(j\times
j)^*\cal B_{ml}$ concentrated on this component, i.e. $$\cal
B_{ml,U}=(j\times j)^*\cal B_{ml}|_{X\times Y\times U}.$$ The above
properties of $\cal B_{ml,U}$ are clear from this definition.

Similarly, we construct the sheaf $\cal B'_{ml,U}$ on $\PP_n(\cal A)\times
\PP_n(\Phi\cal A^{op})$.
It remains to compute the relevant $\circ$-products. This is easily done
using isomorphisms (\ref{iso1}), (\ref{iso2}). Namely, one should start by
computing
$(j\times\id)_*(\cal B_{ml,U}\circ_{\PP_n(\cal A)_U}\cal B'_{ml,U})$ on
$\PP_n(\Phi\cal A)\times\PP_n(\Phi\cal A^{op})_U$. We have \begin{align*}
&(j\times\id)_*(\cal B_{ml,U}\circ_{\PP_n(\cal A)_U}\cal B'_{ml,U})\simeq
((j\times\id)_*\cal B_{ml,U})\circ_{\PP_n(\cal A)_U}\cal B'_{ml,U}\simeq
((\id\times j)^*\cal B_{ml})\circ_{\PP_n(\cal A)_U}\cal B'_{ml,U}\simeq\\
&\cal B_{ml}\circ_{\PP_n(\cal A)}((j\times\id)_*\cal B'_{ml,U})\simeq \cal
B_{ml}\circ_{\PP_n(\cal A)}((\id\times j)^*\cal B'_{ml})\simeq (\id\times
j)^*(\cal B_{ml}\circ_{\PP_n(\cal A)}\cal B'_{ml})\simeq\\ &(\id\times
j)^*(\delta_{\PP_n(\Phi\cal A)})\simeq (j\times\id)_*
(\delta_{\PP_n(\Phi\cal A)_U} )
\end{align*}
Now the situation looks locally as follows: we have an \'etale extension of
NC-algebras $A\ra A_1$, an $A_1\otimes A_1^{op}$-module $M$, and an
isomorphism of $A\otimes A_1^{op}$-modules $M\simeq A_1$. Furthermore, we
have a 2-sided ideal $I\sub A$ such that $A/I$ is commutative and $IM=MI$,
and we know that the induced isomorphism $M/IM\simeq A_1/IA_1$ is an
isomorphism of $A_1/IA_1\otimes A_1/IA_1$-modules (notice that $A_1/IA_1$
is commutative). We claim that in such a situation the above isomoprhism
commutes with the left action of $A_1$. Indeed, the left action of $A_1$ on
$M$ induces a homomorphism $\phi:A_1\ra A_1$ such that $\phi|_A=\id$ and
$\phi\mod IA_1$ is the identity on $A_1/IA_1$. Now the formal \'etaleness
implies that $\phi=\id$.

Thus, we conclude that
$\cal B_{ml,U}\circ_{\PP_n(\cal A)_U}\cal B'_{ml,U}\simeq
\delta_{\PP_n(\Phi\cal A)_U}$ as required. \end{pf}

\subsection{} The sheaf of rings $\O_{\PP_n}$ on $\PP_n$ can be naturally
enlarged as follows.
The central element $t\in\gr_{(n)}(\cal A)_1$ induces a sequence of
embeddings of $\O_{\PP_n}$-bimodules
$$\O_{\PP_n}\ra\O_{\PP_n}(1)\ra\O_{\PP_n}(2)\ra\ldots$$ Now using the
natural morphisms
$$\O_{\PP_n}(m)\otimes_{\O_{\PP_n}}\O_{\PP_n}(l)\ra \O_{\PP_n}(m+l)$$
we can define the ring structure on
the direct limit
$$\wt{\O}_{\PP_n(\cal A)}=
\operatorname{inj.}\lim(\O\ra\O(1)\ra\O(2)\ra\ldots)$$ For example,
if $Y$ is smooth and
$\cal A=D_Y$ is the sheaf of differential operators on $Y$ then
$\wt{\O}_{\PP_{\infty}}$ is the sheaf of (formal) pseudo-differential
operators (the underlying topological space of $\PP_{\infty}$ is the
projectivized cotangent bundle of $Y$). The subsheaf $\O_{\PP_{\infty}}$
consists of (formal) pseudo-differential operators of negative order.

Now let $\PP_1(\cal A)=Y\times Z$ and $U\ra Z$ be an \'etale morphism. Then
one can define invertible $\O_{\PP_n(\cal A)_U}$-bimodules $\O_{\PP_n(\cal
A)_U}(m)$ as follows.
We can regard $\O_{\PP_n(\cal A)}$ as a sheaf on $\PP_n(\cal A)\times
\PP_n(\cal A^{op})$ supported on the diagonal. Let $V$ be a thickening of
the diagonal in $\PP_n(\cal A)\times\PP_n(\cal A^{op})$ on which
$\O_{\PP_n(\cal A)}(m)$ lives. Then there is a canonical \'etale morphism
$V_U\ra V$ and an embedding $V_U\ra \PP_n(\cal A)_U\times\PP_n(\cal
A^{op})_U$. Now by definition $\O_{\PP_n(\cal A)_U}(m)$ is obtained from
$\O_{\PP_n(\cal A)}(m)$ by first pulling back to $V_U$ and then pushing
forward to $\PP_n(\cal A)_U\times\PP_n(\cal A^{op})_U$. It is easy to see
that we still have morphisms of bimodules $\O(m)\ra\O(m+1)$ and
$\O(n)\otimes\O(m)\ra\O(n+m)$ so we can define the algebra
$\wt{\O}_{\PP_n(\cal A)_U}$.

\begin{Thm} In the preceding two theorems one can replace the categories of
$\O$-modules by the categories of $\wt{\O}$-modules.
\end{Thm}

The proof is an application of the analogue of Theorem \ref{main1} for
$D$-modules on NC-schemes.

\section{Noncommutative deformation of the Poincar\'e line
bundle}\label{geom-subsec}

Consider
the following data:
\newline
\noindent
$W$ is a smooth projective variety over $\C$ of dimension $r$,\newline
\noindent
$D\sub W$ is a reduced irreducible effective divisor,\newline \noindent
$V\sub H^0(D,\O_D(D))$ is an $r$-dimensional subspace, such that the
corresponding rational morphism $\phi:D\ra\PP(V^*)$ is generically finite,
\newline \noindent
$U\sub D$ is an open subset such that $\phi|_U$ is \'etale.

From the exact sequence
$$0\ra \O_W\ra\O_W(D)\ra\O_D(D)\ra 0$$
we get a boundary homomorphism
$$V\ra H^0(D,\O_D(D))\ra H^1(W,\O_W).$$
Now let $X$ be the Albanese variety of $W$, $a:W\ra X$ be the Abel-Jacobi map
(associated with some point of $W$). Then we have the canonical isomorphism
$$H^1(X,\O_X)\wt{\ra}
H^1(W,\O_W),$$ in particular, we get a homomorphism
$V\ra H^1(W,\O_W)$. Let
$$0\ra\O_X\ra{\cal E} \ra H^1(X,\O_X)\otimes_{\C}\O_X\ra 0$$ be the
universal extension.
Taking the pull-back of this extension under the map $V\ra H^1(W,\O_W)$ we
obtain an extension
$$0\ra\O_X\ra{\cal E}_V\ra V\otimes_{\C}\O_X\ra 0.$$ Now we define a
commutative sheaf of
algebras on $X$ as follows
$$\cal A_V=\Sym ({\cal E}_V)/(1_{\cal E}-1)$$ where $1_{\cal E}$ is the
image of
$1\in\O_X\ra{\cal E}_V$. Note that $\cal A_V$ is equipped with the
filtration satisfying the
conditions of the previous section. Also by construction we have a
canonical morphism of sheaves
$${\cal E}_V\ra a_*\O_W(D)$$ which induce the homomorphism of
$\O_X$-algebras $${\cal A}_V\ra
a_*\O_W(*D)$$ compatible with natural filtrations, where $\O_W(*D)=
\operatorname{inj.}\lim
\O_W(nD)$.

If the map $V\ra H^1(V,\O_V)$ is an embedding then the dual $D$-algebra
$\Phi{\cal A}_V$
is the algebra of differential operators ``in directions $V$", where we
consider $V$ as a subspace
in $H^1(X,\O_X)\simeq H^0(\hat{X},T_{\hat{X}})$.
By Theorem \ref{main1} we get a functor from the derived category of ${\cal
A}_V$-modules to
the derived category of $\Phi{\cal A}_V$-modules. We can restrict this
functor to the category of
$\O_W(*D)$-modules. For example, if $D$ is ample the Fourier transform of
$\O_W(*D)$ is a coherent
$\Phi{\cal A}_V$-module. In the case when $W$ is a curve, $D$ is a point,
and $a(D)=0\in X$ one can
show that the latter $\Phi{\cal A}_V$-module is free of rank 1 at general
point. However, in
general $\cal F(\O_W(*D))$ is not free as $\Phi\cal A_V$-module even at
general point unless
$a(D)=0$. To get a module which is free of rank 1 at general point we have
to pass to
microlocalization and use \'etale localization "in vector fields direction"
as described below.

Let $D_{(n)}\sub W$ be the closed subscheme corresponding to the divisor
$nD$. Then
$\Proj(\gr_{(n)}(\O_W(*D)))\simeq D_{(n)}$, hence by functoriality we have
a morphism
$$a_n:D_{(n)}\ra\PP_n(\cal A_V)$$
and an isomorphism
$$a_{n,*}\O_{D_{(n)}}\simeq\ml_n(\O_W(*D))$$ where $\O_W(*D)$ is considered
as a
$\cal A_V$-module.

Let us start with the case $n=1$. Note that
$$a_1:D\ra \PP_1(\cal A_V)\simeq X\times \PP(V^*)$$
is the natural map induces by $a$ and by
$\phi$. Hence, applying Fourier-Mukai transform to $a_{1,*}\O_D$ over a
general point
of $\PP(V^*)$
one gets a free module of rank equal to the degree of $\phi$. To get a free
module of rank 1
at
general point we use the \'etale base change $U\ra\PP(V^*)$. Namely, we
replace $\PP_1(\cal
A_V)=X\times\PP(V^*)$ by
$\PP_1(\cal A_V)_U=X\times U$ and $a_1$ by the morphism
$$a^U_1:U\ra\PP_1(\cal A_V)_U=X\times U:u\mapsto (a(u),u).$$
Then the following lemma is clear.

\begin{Lem}   The  relative
Fourier transform of $a^U_{1,*}\O_U$ is the line bundle $(id\times
a|_U)^*({\cal P})$ on
$\hat{X}\times U$, where $\cal P$ is the Poincar\'e line bundle. \end{Lem}

Now let $\PP_n(\cal A_V)_U$ be the \'etale scheme over $\PP_n(\cal A_V)$
which is a thickening of $X\times U$. Let also $U_{(n)}$ be the open subset
of $D_{(n)}$ which is a nilpotent thickening of $U$.
Then we have a commutative diagram
\begin{equation}
\begin{array}{ccc}
U & \lrar{} & \PP_n(\cal A_V)_U \\
\ldar{} & & \ldar{} \\
U_{(n)} & \lrar{} & \PP_n(\cal A_V)
\end{array}
\end{equation}
where the top horizontal arrow is the composition of $a_1^U$ and the closed
embedding
$\PP_1(\cal A_V)_U\ra\PP_n(\cal A_V)_U$, the bottom horizontal arrow is the
restriction of $a_n$. Since the right vertical morphism is \'etale this
diagram gives rise to a morphism
$$a^U_n: U_{(n)}\ra\PP_n(\cal A_V)_U$$
filling the diagonal in the above commutative square. It is easy to check
that $a_n^U|_{U_{(n-1)}}=a_{n-1}^U$. Now we define the sequence of coherent
sheaves on $\PP_n(\Phi\cal A_V)_U$ by setting
$$\cal L_n=\Phi_{(n)}(a^U_{n,*}\O_{U_{(n)}}).$$ These sheaves satisfy
$\cal L_{n+1}|_{\PP_n}\simeq \cal L_n$, hence we can consider the
projective limit $\cal L_{\infty}$ of $\cal L_n$ which is a coherent sheaf
on the formal NC-scheme
$\PP_{\infty}(\Phi\cal A_V)_U$.

\begin{Prop} The $\O_{\PP_{\infty}(\Phi\cal A_V)_U}$-module $\cal
L_{\infty}$ is locally free of rank 1. \end{Prop}
\begin{pf} One has an exact sequence
$$0\ra a_{n-1,*}^U\cal O_{U_{(n-1)}}(-1)\ra a_{n,*}^U\cal O_{U_{(n)}}\ra
a_{1,*}^U\cal O_U\ra 0.$$ Applying the functor $\Phi_{(n)}$ and passing to
the limit we obtain an exact sequence
$$0\ra \cal L_{\infty}(-1) \stackrel{t}{\ra} \cal L_{\infty}\ra \cal L_1\ra 0$$
where $t$ is the canonical central element in $\cal O(1)$. It remains to
use the following simple algebraic fact. Let $A$ be a noetherian ring,
$t\in A$ be a non zero divisor such that $At=tA$ and 
$A=\operatorname{proj.}\lim A/t^nA$.
Let $M$ be a finitely generated left $A$-module such that $t$ is not a
divisor of zero in $M$ and $M/tM$ is a free $A/tA$-module
of rank 1. Then $M$ is a free $A$-module of rank 1. \end{pf}

Notice that in the case when $W=C$ is a curve $D=P$ is a point, $U=D$ the
module $\cal L_{\infty}$ was used in [R1] to construct Krichever's solution
to the KP hierarchy. In this case $\cal L_{\infty}$ is a locally free
module of rank-1 over the microlocalization of
the subalgebra $\O[\xi]$ in the ring of the differential operators on the
Jacobian $J(C)$ generated by the vector field $\xi$ which comes from the
boundary homomorphism $H^0(O_P(P))\ra H^1(\O_C)$. The key point is that
$\cal L_{\infty}$ has also an action of completion of $\O_C(*P)$ at $P$
(which is isomorphic to the ring of Laurent series) commuting with the
action of pseudo-differential operators in $\xi$.

\bigskip

\noindent
{\bf References}

\bigskip

\noindent
[BB] A.~Beilinson, J.~Bernstein, {\it A proof of
Jantzen conjectures}. I.~M.~Gelfand Seminar, 1--50, Adv. Soviet Math., 16,
Part 1,
Amer. Math. Soc., Providence, RI, 1993.

\noindent
[BD] A.~Beilinson, V.~Drinfeld, {\it Quantization of
Hitchin's fibration and Langlands' program}. Algebraic and geometric
methods in mathematical physics (Kaciveli, 1993), 3--7, Math. Phys. Stud.,
19, Kluwer Acad. Publ., Dordrecht, 1996.

\noindent
[BC] J.~Burchnall, T.~Chaundy, {\it Commutative ordinary differential
operators},
 Proc. London Math. Soc., 21 (1923), 420--440;
{\it Commutative ordinary differential
operators II}, Proc. Royal Soc.
London (A), 134 (1932), 471--485.

\noindent
[Kr] I.M.~Krichever,
{\it Algebro-geometric construction of
the Zaharov-Shabat equations and their periodic solutions}.
Soviet Math. Dokl. 17 (1976), 394--397;
{\it Integration of nonlinear equations by the
methods of nonlinear geometry}, Funk. Anal. i Pril, 11  (1977),  15-- 31.

\noindent
[L] G.~Laumon, {\it Transformation de Fourier generalisee}, preprint
alg-geom 9603004.

\noindent
[K] M.~Kapranov, {\it Noncommutative geometry based on commutator
expansions}, preprint math.AG/9802041.

\noindent
[M] S.~Mukai, {\it Duality between $D(X)$ and $D(\hat{X})$ with its
application to Picard sheaves}.
 Nagoya Math. J. 81 (1981), 153--175.

\noindent
[R1] M.~Rothstein, {\it Connections on the Total Picard Sheaf and the KP
Hierarchy}, Acta Applicandae Math. 42 (1996), 297--308.

\noindent
[R2] M.~Rothstein, {\it Sheaves with connection on abelian varieties}, Duke
Math. Journal 84 (1996), 565--598.

\vspace{3mm}

{\sc Department of Mathematics, Harvard University, Cambridge,
 MA 02138

Department of Mathematics, University of Georgia, Athens, GA 30602}

{\it E-mail addresses:} apolish@@math.harvard.edu, rothstei@@math.uga.edu

\end{document}